\documentclass[11pt]{article}

\usepackage[T1]{fontenc}
\usepackage{amsmath,amssymb,amsthm}
\usepackage{mathtools}
\usepackage{geometry}
\usepackage{graphicx}
\graphicspath{{figs/}}
\usepackage{booktabs}
\usepackage{array}
\usepackage{longtable}
\usepackage{enumitem}
\usepackage{xcolor}
\usepackage[colorlinks=true,linkcolor=blue!50!black,citecolor=blue!50!black,urlcolor=blue!50!black]{hyperref}

\newtheorem{theorem}{Theorem}[section]
\newtheorem{lemma}[theorem]{Lemma}
\newtheorem{proposition}[theorem]{Proposition}
\newtheorem{corollary}[theorem]{Corollary}

\theoremstyle{remark}
\newtheorem{remark}[theorem]{Remark}
\newtheorem*{proof*}{Proof}

\newcommand{\R}{\mathbb{R}}

\newcommand{\Z}{\mathbb{Z}}
\newcommand{\tr}{\operatorname{tr}}
\newcommand{\diag}{\mathrm{diag}}

\newcommand{\Spec}{\operatorname{Spec}}
\newcommand{\ind}{\operatorname{ind}}
\newcommand{\St}{\operatorname{St}}
\newcommand{\Gr}{\operatorname{Gr}}
\newcommand{\dd}{\mathrm{d}}

\title{\textbf{Hessian--Maslov Geometry of Oja--Brockett Flows}}

\author{Shintaro Yoshizawa\\[2pt]
\normalsize Nagoya Mathematical and Information Science Research\\
\normalsize \href{mailto:shintaro.yoshizawa.net@gmail.com}{shintaro.yoshizawa.net@gmail.com}}
\date{August 13, 2026}

\begin{document}
\maketitle

\begin{abstract}
This paper is organized in three parts. Part~I gives a careful account of the
theoretical foundations of the \textbf{Maslov index} associated with degenerate
points (caustics) of the Legendre map --- its origin in the breakdown of phase
in the semiclassical approximation, its geometric formulation as an
intersection number in the Lagrangian Grassmannian, H\"ormander's local
formula, and the dynamical viewpoint of bifurcating fixed points of a gradient
flow --- and then works out in detail the simplest example, the
two-dimensional vector-valued quartic potential
$f(\xi)=-\tfrac12\xi^\top A^2\xi+\tfrac14(\xi^\top A\xi)^2$ ($A$: a $2\times2$
positive-definite symmetric matrix). We determine its critical-point structure
and its caustic (a cusp at $\lambda_2=3\lambda_1$), and we exhibit the
bifurcation in which minima and saddles exchange roles as the eigenvalue
ratio crosses $1$.

Part~II lifts this theory in earnest from $\xi\in\R^2$ to a tall $n\times k$
real matrix $X$. Introducing the potential
$f(X)=-\tfrac12\tr(A^2XB^2X^\top)+\tfrac14\tr[(AXBX^\top)^2]$, built from a
positive-definite symmetric matrix $A$ and a diagonal positive-definite matrix
$B=\diag(b_1,\dots,b_k)$, we confirm numerically that at the fixed points of
its negative gradient flow the diagonal entries of $X^\top AX$ align with the
top $k$ eigenvalues of $A$, exactly as predicted by \textbf{Brockett's sorting
theorem} (1991). When $A$ or $B$ has degenerate eigenvalues we show that the
fixed-point set expands from isolated points into a \textbf{flag manifold},
and we confirm numerically that when an eigenvalue of $A$ crosses a wall of
the Weyl chamber (the locus where two eigenvalues coincide) the Hessian
degenerates and the sort order exchanges --- the exact generalization of the
stability exchange at $\lambda_1=\lambda_2$ seen in Part~I. We situate this
wall-crossing phenomenon within the intersection theory of Lagrangian
submanifolds and the Maslov cycle, and show that the combinatorial inversion
number of permutations and the length function of the Coxeter group, on the
one hand, and the differential-geometric Maslov index, on the other, are two
faces of the same phenomenon.

Part~III is the paper's principal contribution. We give a \textbf{rigorous}
analysis of the structure underlying the phenomena observed numerically in
Part~II --- Brockett sorting, flag manifolds, and Hessian degeneracy at Weyl
walls --- by computing the ambient Hessian of the general Oja--Brockett
gradient flow $\dot X=AXB-XBX^\top AX$ associated with two independent SPD
matrices $A,B$. We decompose the Hessian completely into normal, vertical,
and horizontal sectors and give closed spectral formulas depending on both
$A$ and $B$. In particular, the ``identically vanishing gradient along the
flag-manifold direction,'' which Part~II observed numerically only for
degeneracies of $A$, turns out to be a special case of the general formula
$\mu^V_{ij}=(\lambda_i-\lambda_j)(b_i-b_j)$; we show that degeneracy of $B$
($b_i=b_j$) independently triggers the very same phenomenon (a \emph{$B$-wall}).
That is, we uncover, within a structure that Part~II could only glimpse
implicitly, a \textbf{second wall system} dual to the Weyl wall on the $A$
side, and we establish, as a theorem, the mechanism by which a genuine
crossing of eigenvalues becomes an avoided crossing under this double-wall
structure, together with fully explicit closed-form coupling constants
(horizontal: $\gamma=-b_jC_{qr}$; vertical:
$\gamma^V=(\lambda_m^2-\lambda_i\lambda_l)/\sqrt{(\lambda_i+\lambda_m)(\lambda_m+\lambda_l)}$).
\end{abstract}

\tableofcontents

\clearpage
\section{Introduction}

\subsection{Motivation and Problem Statement}

Two classical bodies of mathematics rarely appear side by side. The first is
the theory of \textbf{caustics and the Maslov index}, born from the
semiclassical (WKB) approximation of quantum mechanics: the graph of a
gradient map is always a Lagrangian submanifold of phase space, and the
points at which this graph fails to project bijectively onto configuration
space --- the caustics --- are counted, with sign, by an integer invariant,
the Maslov index. The second is the theory of \textbf{matrix gradient
flows} from control theory and numerical linear algebra: Brockett's
double-bracket flow sorts the eigenvectors of a symmetric matrix along a
Stiefel manifold, converging almost surely to the eigenbasis ordered by a
weight matrix. On the surface these look like unrelated subjects --- one
analytic and geometric, concerned with the propagation of a wave; the other
algebraic and dynamical, concerned with an algorithm for diagonalization.

The present paper's point of departure is the observation that these two
theories are, in a precise sense, the same phenomenon viewed at different
dimensions. The Hessian of a potential is symmetric; the graph of its
gradient is therefore always Lagrangian, whether the potential lives on
$\R^2$ or on a space of $n\times k$ matrices. A degenerate point of the
Hessian --- a caustic, in the language of Part~I --- is exactly a fixed
point at which a matrix gradient flow loses hyperbolicity, exchanges the
stability of two critical points, and reorders the eigenvectors it sorts.
This identification lets us import the full analytic machinery of the
Maslov index (Lagrangian Grassmannians, H\"ormander's local formula,
intersection numbers with the Maslov cycle) into the study of Brockett-type
sorting flows, and, conversely, lets the combinatorics of the Weyl group and
of flag manifolds illuminate what a ``degenerate point'' means once the
scalar potential $f(\xi)$ is replaced by a genuinely matrix-valued one.

The specific dynamical system that carries this identification is the
\textbf{Oja--Brockett gradient flow}
\[
\dot X=AXB-XBX^\top AX,\qquad X\in\R^{n\times k},
\]
associated with two independent symmetric positive-definite matrices $A$
(an $n\times n$ ``data'' matrix) and $B$ (a $k\times k$ ``weight'' matrix).
The special case $B=I$ is classical: it is Oja's principal-subspace flow,
and its fixed points, sorted by Brockett's theorem, are governed entirely by
the eigenvalue structure of $A$. But $B=I$ collapses one of the two
genuinely independent degrees of freedom of the flow, and the existing
literature --- including the numerical study we ourselves carry out in
Part~II of this paper --- had never asked what happens to the Hessian, and
hence to the associated Maslov geometry, once $B$ is allowed to be a
general diagonal SPD matrix with distinct entries. This is the central
question the paper answers: \emph{does the weight matrix $B$ generate its
own, independent mechanism of Hessian degeneracy, on a par with the
classical Weyl-chamber wall generated by $A$; and if so, what governs the
transition, under either mechanism, from a genuine crossing of eigenvalues
to an avoided crossing?}

\subsection{Summary of Contributions}

The paper answers this question completely, and organizes the answer as a
progression from the simplest possible instance of the phenomenon to its
fully general matrix form. The main original contributions are:

\begin{enumerate}[label=(\roman*)]
\item \textbf{A complete classification of the minimal (scalar/vector)
model.} For the two-dimensional quartic potential
$f(\xi)=-\tfrac12\xi^\top A^2\xi+\tfrac14(\xi^\top A\xi)^2$, we give closed
forms for its caustic curve, its cusp locus ($\lambda_2=3\lambda_1$), and
the associated Maslov index, and we embed $f$ within two natural
one-parameter deformation families, $h_\alpha$ (Box--Cox/Tsallis type) and
$\Phi_\alpha$ (pure-power type), showing that only the latter contains $f$
as an exact special case ($\Phi_2\equiv f$) and correctly inherits its
scale invariance (Sections~3--6).

\item \textbf{A numerical atlas of the matrix-valued phenomenon.} We verify,
for the $n\times k$ generalization of $f$, that Brockett's sorting theorem,
the emergence of flag-manifold degeneracy, and Hessian degeneracy at
Weyl-chamber wall crossings all persist, and we identify the precise sense
in which wall-crossing is a higher-dimensional Maslov-cycle intersection
(Sections~8--14). This numerical study is deliberately restricted to
$B=I$ and is designed to isolate exactly the questions that Part~III must
answer rigorously.

\item \textbf{The rigorous Hessian--Maslov theory of the general $(A,B)$
flow (the paper's principal contribution).} We decompose the ambient
Hessian of $V_{A,B}$ at every Stiefel-manifold equilibrium into three
mutually non-interacting sectors --- normal, vertical, horizontal --- and
give closed-form spectra for each, valid for arbitrary SPD $A,B$
(Theorem~20.1). This decomposition reveals a structure invisible at $B=I$:
the vertical sector carries a \emph{second} wall system,
$\{b_i=b_j\}$ (the \emph{$B$-wall}), governed by the closed formula
$\mu^V_{ij}=(\lambda_i-\lambda_j)(b_i-b_j)$ (Lemma~17.1), independent of and
coupled to the classical Weyl wall $\{\lambda_j=a_q\}$ on the $A$ side
(Sections~15--21).

\item \textbf{Explicit avoided-crossing theorems for both wall systems.}
We prove that a genuine eigenvalue crossing becomes an avoided crossing
under an explicit, fully closed-form coupling constant in each sector ---
$\gamma=-b_jC_{qr}$ on the horizontal (Weyl-wall) side (Theorem~26.1) and
$\gamma^V=(\lambda_m^2-\lambda_i\lambda_l)/\sqrt{(\lambda_i+\lambda_m)(\lambda_m+\lambda_l)}$
on the vertical ($B$-wall) side (Corollary~26.4) --- and we show that the
two mechanisms are structural duals of one another, coupling through a
shared selected column ($j$, respectively $m$) rather than through a shared
non-selected or selected index pair (Sections~22--26).
\end{enumerate}

\subsection{Organization of the Paper}

The paper is organized in three parts, of increasing generality, each
building the mathematical apparatus needed for the next.

\textbf{Part~I} (Sections~2--7) develops the Maslov index in its classical,
scalar setting and works the minimal vector-valued example by hand.
Section~2 recalls the origin of the Maslov index in the WKB approximation,
its formulation as an intersection number in the Lagrangian Grassmannian,
H\"ormander's local formula, and its reinterpretation as a bifurcation
invariant for gradient-flow fixed points. Section~3 computes, for the
two-dimensional potential $f(\xi)$, its caustic curve, its cusp locus, and
its Maslov index along explicit paths, and exhibits the elementary
bifurcation --- the exchange of minima and saddles as $\lambda_2/\lambda_1$
crosses $1$ --- that is the seed of everything that follows. Sections~4--7
situate $f$ within the two deformation families $h_\alpha$ and $\Phi_\alpha$,
establish that $f=\Phi_2$ exactly, and contrast the two families'
qualitatively different degeneracy mechanisms.

\textbf{Part~II} (Sections~8--14) lifts the vector $\xi\in\R^2$ to an
$n\times k$ matrix $X$ and numerically maps out the resulting phenomena
before they are proved. Sections~8--9 set up four natural matrix-valued
potential families and their gradients. Section~10 confirms Brockett's
sorting theorem numerically for the flow generated by $f(X)$. Section~11
shows that eigenvalue degeneracy of either $A$ or $B$ expands isolated
fixed points into flag manifolds. Section~12 confirms Hessian degeneracy
and sort-order exchange at Weyl-chamber wall crossings and interprets this
as a Maslov-cycle intersection in the matrix setting. Section~13 uncovers,
in a related penalized potential $\varphi_\alpha$, a sign-driven bifurcation
between principal- and minor-component sorting. Section~14 catalogues the
combinatorics of critical points for general $(n,k)$.

\textbf{Part~III} (Sections~15--29) is the paper's theoretical core, and
supplies the rigorous proofs that Part~II's numerics call for. Section~15
sets up the general $(A,B)$ Oja--Brockett potential, its gradient flow, and
the threefold normal/vertical/horizontal decomposition of the tangent space
at a Stiefel-manifold equilibrium. Sections~16--19 compute the Hessian on
each sector and on their mutual couplings, culminating in the complete
spectral formula of Section~20 and the double-wall structure articulated in
Section~21. Sections~22--26 develop the crossing-to-avoided-crossing
transition theorem for both the horizontal (Weyl-wall) and vertical
($B$-wall) sectors, with fully explicit coupling constants. Section~27
revisits the minimal model of Part~I as the degenerate case $k=1$ of this
general theory, Section~28 records open problems, and Section~29 concludes.

\part{Foundations of the Maslov Index and the Minimal Vector Model}

Building the general theory of Part~III requires, first, a firm grasp of
the classical, scalar Maslov index and of the one nontrivial phenomenon ---
a stability-exchange bifurcation --- that a two-dimensional vector-valued
potential can already exhibit. This part supplies both. Section~2 recalls
the origin, geometric formulation, and computational content of the Maslov
index; Section~3 works through the minimal example $f(\xi)$ in complete
detail, including its Legendre-transformed front surface and an explicit
Maslov-index computation along a radial path; Sections~4--7 place $f$ within
two broader deformation families in order to isolate exactly which of its
features are structural (shared by the whole family) and which are
accidents of the particular normalization $f=\Phi_2$. The matrix
generalization of Part~II, and the rigorous theory of Part~III, both trace
every one of their central phenomena back to a construction introduced here.

\section{What Is the Maslov Index? Definition, Meaning, and Applications}

\subsection{Origin: the breakdown of phase in the semiclassical (WKB) approximation}

The Maslov index is a concept introduced by V.\,P.\ Maslov in 1965 in the
course of rigorously formulating the semiclassical (WKB) approximation of
quantum mechanics (with the pioneering 1958 work of J.\,B.\ Keller as a
precursor). One approximates the solution of the Schr\"odinger equation, in
terms of the classical action $S(x)$, by
\[
\psi(x)\sim\frac{C}{\sqrt{|\det(\partial^2S/\partial x^2)|}}\exp\!\Big(\frac{i}{\hbar}S(x)\Big)
\]
(the WKB approximation). Here $S(x)$ is the projection onto $x$-space of the
Lagrangian submanifold $\Lambda=\{(x,\partial S/\partial x)\}$ (the collection
of classical trajectories) inside the phase space $T^*\R^n$.

The difficulty is that this projection need not be invertible everywhere.
At points where $\Lambda$ fails to project smoothly and bijectively onto
$x$-space (points we call ``degenerate points'' or ``caustics'' in this
paper), the Hessian $\partial^2S/\partial x^2$ becomes singular and the
amplitude $1/\sqrt{|\det(\cdots)|}$ diverges. Maslov gave a procedure for
skirting such a caustic and correctly gluing the wave function back together,
showing that the phase of the wave function picks up an extra shift of
$-\pi/2$ at each passage through a caustic. The integer obtained by dividing
the accumulated total of this phase shift by $\pi/2$ is the Maslov index:

\begin{quote}
\textbf{An integer that counts the ``extra'' quantum-mechanical phase shift
that a semiclassical wave function must accumulate, each time it passes a
caustic, in order to remain well defined along a path.}
\end{quote}

\begin{figure}[htbp]
\centering
\includegraphics[width=0.92\textwidth]{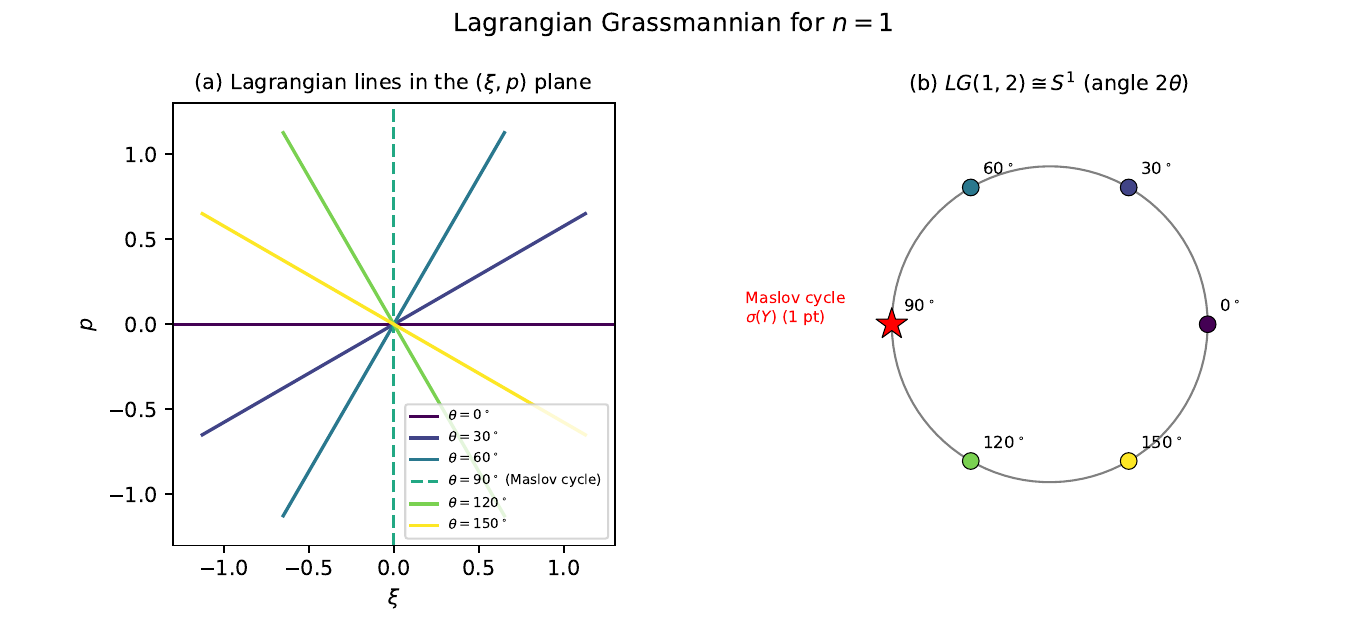}
\caption{The Lagrangian Grassmannian for $n=1$. (a) Various lines through the
origin (all Lagrangian) in the $(\xi,p)$ plane; (b) each line represented as a
point on the circle $S^1$. The vertical line ($\theta=90^\circ$, red star)
corresponds to the Maslov cycle.}
\end{figure}

\subsection{Geometric formulation: intersection number in the Lagrangian Grassmannian}

Treating the gradient $p=\nabla V(\xi)$ of a function $V(\xi_1,\xi_2)$ of two
variables as an independent variable, the four-dimensional space
$(\xi_1,\xi_2,p_1,p_2)=T^*\R^2$ (phase space) carries the symplectic form
\[
\omega=\dd\xi_1\wedge \dd p_1+\dd\xi_2\wedge \dd p_2.
\]
A two-dimensional subspace $V\subset T^*\R^2$ with $\omega|_V\equiv0$ is
called a \textbf{Lagrangian subspace}.

\begin{proposition}[Lagrangian property of the graph of a gradient map; Ekeland (1977)]
Let $V:\R^n\to\R$ be $C^\infty$. Then the graph of the gradient map,
$\Lambda_V=\{(\xi,\nabla V(\xi)):\xi\in\R^n\}\subset T^*\R^n$, is always a
Lagrangian submanifold.
\end{proposition}

\begin{proof*}
For $u=(\dd\xi,H\,\dd\xi)$, $v=(\delta\xi,H\,\delta\xi)$ (where
$H=\partial^2V/\partial\xi^2$),
\[
\omega(u,v)=\dd\xi^\top H\delta\xi-(H\,\dd\xi)^\top\delta\xi
=\dd\xi^\top H\delta\xi-\dd\xi^\top H^\top\delta\xi=0,
\]
which vanishes identically simply because $H$ is symmetric (equality of
mixed partials, $\partial^2V/\partial\xi_i\partial\xi_j
=\partial^2V/\partial\xi_j\partial\xi_i$). \qed
\end{proof*}

\begin{quote}
From nothing more than the elementary calculus fact that the Hessian is a
symmetric matrix, the graph $\Lambda_V$ of a gradient map is automatically a
Lagrangian submanifold.
\end{quote}

This is the fact we repeatedly invoke throughout the paper as ``the graph of
a gradient map is Lagrangian.'' Ekeland proves this rigorously, without any
appeal to convexity, using the one-form $\omega_z=\dd z-\sum_i p_i\,\dd\xi_i$
on $\R^n\times\R^n\times\R$ (this is the theoretical underpinning of the
claim, used throughout the paper, that the Lagrangian character of
$\Lambda_V$ itself survives even near degenerate points).

The simplest example ($n=1$): in the phase space $\R^2=(\xi,p)$, every line
through the origin is a Lagrangian subspace (tangent vectors along a
one-dimensional subspace are automatically parallel, so the signed area
automatically vanishes). Hence the space of all Lagrangian subspaces,
$LG(1,2)$, can be identified with the set of all such lines, parametrized by
an angle $\theta\in[0,\pi)$ --- that is, with the circle $S^1$ (represented
via the doubled angle $2\theta$). Among these, only the vertical line
($\theta=90^\circ$) can never be represented as the graph of any
single-variable function; we call this the \textbf{Maslov cycle}. As
$V''(\xi_0)\to\pm\infty$ the graph approaches this vertical line, coinciding
with it exactly when $V''=0$ (a caustic).

For general $n$ (in this paper, mainly $n=2$, and its generalization to
$n\times k$ matrices), the Maslov cycle is a codimension-one submanifold of
$LG(n,2n)$, and the striking fact (Arnold) that the fundamental group is
\[
\pi_1(LG(n,2n))\cong\Z
\]
holds independently of dimension. The signed count of intersections of the
family of Lagrangian planes swept out along a path parametrized by $s$ with
the Maslov cycle is the Maslov index. H\"ormander's local formula
\begin{equation}
\ind_{s_0}=\frac{\operatorname{sgn}\det H(s_0^+)-\operatorname{sgn}\det H(s_0^-)}{2}
\end{equation}
translates this ``signed intersection number'' into the computable language
of the sign change of a Hessian determinant.

\subsection{What does the Maslov index tell us?}

\begin{enumerate}[label=(\Alph*)]
\item How many times the path has passed through a degenerate point --- an
instant at which the gradient map loses local invertibility (a bifurcation
point where several $\xi$ give the same slope $p$).
\item Whether that passage is a ``round trip'' or a ``one-way'' passage ---
a round trip through the same degenerate point cancels the signs and
contributes net zero.
\item Whether the start and end points of the path lie in the same
``chamber'': the complement $LG(n,2n)-\sigma(Y)$ of the Maslov cycle
splits into several connected components, and a nonzero net index means the
start and end points lie in different chambers (the signature of the Hessian
has changed).
\end{enumerate}

\subsection{Another viewpoint: bifurcation of the fixed points of a gradient flow}

For the gradient flow $\dot\xi=-\nabla V(\xi)$ of a potential $V(\xi)$, the
fixed points are the critical points of $V$, and their stability is
determined by the sign of the eigenvalues of the Hessian. Whereas the caustic
condition $\det H(\xi)=0$ is defined at every point of the $\xi$-plane
(without requiring $\nabla V=0$), the bifurcation of fixed points of the
gradient flow occurs only on the more special subset
\[
\{\xi:\det H(\xi)=0\}\supset\{\xi:\det H(\xi)=0\ \text{and}\ \nabla V(\xi)=0\}.
\]
We use this viewpoint in Section~3, in the concrete example of a stability
exchange among the fixed points of an isotropic quartic potential.

\subsection{Concrete applications}

The Maslov index is used throughout physics, mathematics, and engineering:
the \textbf{Bohr--Sommerfeld--Maslov quantization condition} (EBK
quantization, $\oint p\,\dd x=2\pi\hbar(n+\mu/4)$; the origin of the
zero-point energy of the harmonic oscillator), \textbf{Sturm oscillation
theory and the Morse index theorem} (as the number of conjugate points of a
one-dimensional Sturm--Liouville-type equation), \textbf{geometric-optical
caustics} (the actual $\pi/2$ phase jump at rainbows and focal points),
\textbf{catastrophe theory} (a cusp is a self-intersection point of the Maslov
cycle itself), and, finally, as a \textbf{diagnostic for statistical
identifiability}: the Legendre map of an exponential family (natural
parameter $\to$ mean parameter) never touches the Maslov cycle, because its
dispersion function is always positive, so its Maslov index vanishes
identically; but for non-classical \emph{deformed} exponential families such
as $h_\alpha$ or $\Phi_\alpha$, regions where this one-to-one property breaks
down genuinely exist, and the Maslov index serves as a quantitative
diagnostic of the ``location,'' ``multiplicity,'' and ``orientation'' of that
breakdown.

\section{The Minimal Model: Degenerate-Point Structure of the Vector-Valued Potential $f(\xi)$}

We now apply the general theory above to the simplest concrete example ---
one that reappears in Part~III as the case $n=2,k=1,B=1$ of the general
$(A,B)$ matrix theory.

\subsection{The potential and the caustic condition}

Let $A$ be a $2\times2$ positive-definite symmetric matrix and consider
\begin{equation}
f(\xi)=-\frac12\xi^\top A^2\xi+\frac14(\xi^\top A\xi)^2. \label{eq:f-def}
\end{equation}
In the eigenbasis of $A$ (with $A=\diag(\lambda_1,\lambda_2)$), writing
$Q:=\xi^\top A\xi=\lambda_1\xi_1^2+\lambda_2\xi_2^2$, we have
$f(\xi)=-\tfrac12(\lambda_1^2\xi_1^2+\lambda_2^2\xi_2^2)+\tfrac14Q^2$
(Figure~2 shows the isotropic shape).

\begin{figure}[htbp]
\centering
\includegraphics[width=0.5\textwidth]{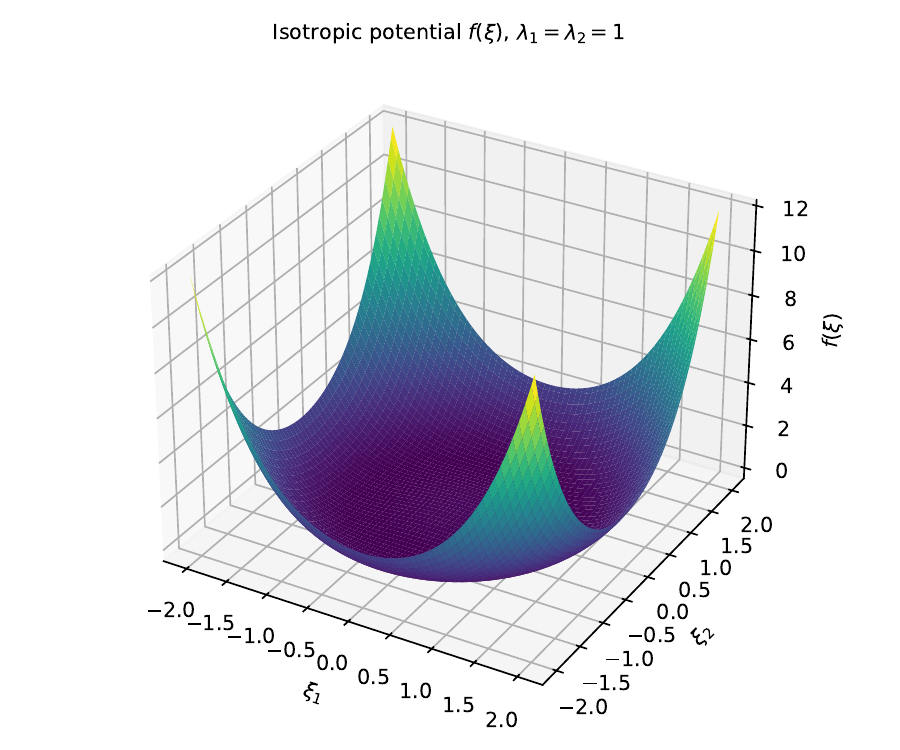}
\caption{Shape of the potential $f(\xi)$ in the isotropic case ($\lambda_1=\lambda_2=1$).}
\end{figure}

The condition for the Hessian determinant to vanish (the caustic condition)
can be written in the coordinate-free, invariant form
\begin{equation}
3(\xi^\top A\xi)^2-\xi^\top A^2\xi-3(\det A)|\xi|^2+\det A=0. \label{eq:caustic-cond}
\end{equation}
In the isotropic case $\lambda_1=\lambda_2=\lambda$, writing $\rho=|\xi|$, this
becomes $3\rho^4-4\rho^2+1=0\iff(3\rho^2-1)(\rho^2-1)=0$, so the caustic
curve, independently of $\lambda$, always decomposes into two universal
(scale-invariant) concentric circles
\[
\rho_T=1\ \text{(outer, valley bottom and degenerate)},\qquad
\rho_R=\frac{1}{\sqrt3}\ \text{(inner, a pure fold)}
\]
(Figure~3).

\begin{figure}[htbp]
\centering
\includegraphics[width=0.45\textwidth]{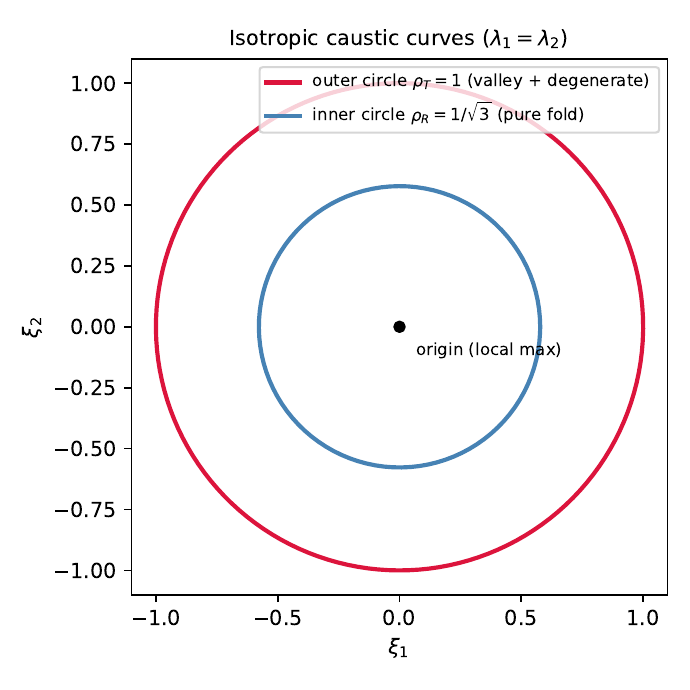}
\caption{Isotropic caustic curves ($\lambda_1=\lambda_2$). Independently of
$\lambda$, the radii are always $\rho_T=1$, $\rho_R=1/\sqrt3$.}
\end{figure}

In the anisotropic case, the caustic deforms into a smooth quartic curve
(Figure~4); precisely when the ratio is $\lambda_2=3\lambda_1$ (a
proportionality passing through the origin) does a genuine cusp form, at
$\xi=(0,\pm1/\sqrt3)$ (Figure~5).

\begin{figure}[htbp]
\centering
\includegraphics[width=\textwidth]{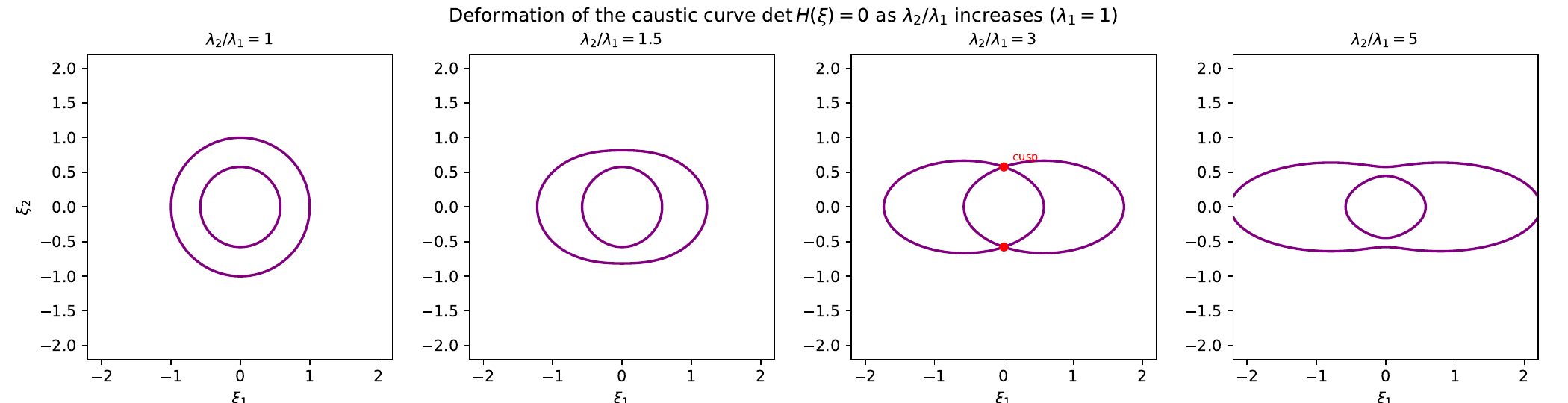}
\caption{Deformation of the caustic curve $\det H(\xi)=0$ as $\lambda_2/\lambda_1$
increases from $1\to1.5\to3\to5$, with $\lambda_1=1$ fixed.}
\end{figure}

\begin{figure}[htbp]
\centering
\includegraphics[width=0.45\textwidth]{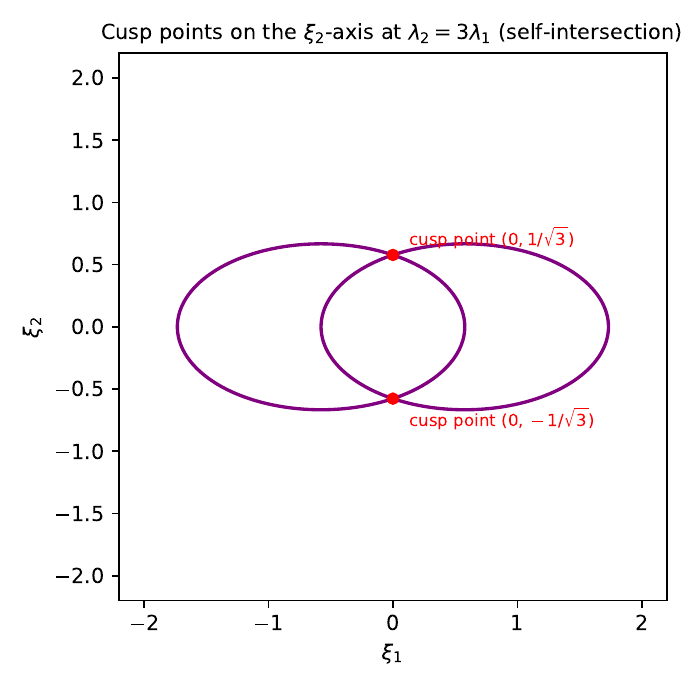}
\caption{Cusp points on the $\xi_2$-axis at $\lambda_2=3\lambda_1$ (a self-intersection).}
\end{figure}

\subsection{Critical-point structure as fixed points of the gradient flow}

Solving $\nabla f(\xi)=0$ (components $p_1=\lambda_1\xi_1(Q-\lambda_1)$,
$p_2=\lambda_2\xi_2(Q-\lambda_2)$), for $\lambda_1\neq\lambda_2$ we find five
critical points $\xi=(0,0),(\pm1,0),(0,\pm1)$. The Hessians are
\[
H(0,0)=\begin{pmatrix}-\lambda_1^2&0\\0&-\lambda_2^2\end{pmatrix}
\ (\text{always a local max}),\qquad
H(\pm1,0)=\begin{pmatrix}2\lambda_1^2&0\\0&\lambda_2(\lambda_1-\lambda_2)\end{pmatrix},
\]
\[
H(0,\pm1)=\begin{pmatrix}\lambda_1(\lambda_2-\lambda_1)&0\\0&2\lambda_2^2\end{pmatrix}.
\]
The second eigenvalue $\lambda_2(\lambda_1-\lambda_2)$ of $H(\pm1,0)$ and
the first eigenvalue $\lambda_1(\lambda_2-\lambda_1)$ of $H(0,\pm1)$ are
exactly opposite in sign. Hence:
\begin{itemize}
\item $\lambda_1>\lambda_2$: $(\pm1,0)$ are stable minima, $(0,\pm1)$ are saddles.
\item $\lambda_1<\lambda_2$: the roles swap --- $(0,\pm1)$ are minima, $(\pm1,0)$ are saddles.
\end{itemize}
As the eigenvalue ratio $\lambda_2/\lambda_1$ crosses $1$, the gradient
flow's two stable minima and two saddles exchange in a bifurcation exactly at
$\lambda_1=\lambda_2$. At $\lambda_1=\lambda_2$ both $H(\pm1,0)$ and
$H(0,\pm1)$ have exactly one eigenvalue vanish (degenerate), and the four
points merge continuously onto the ``valley-bottom circle'' $\rho=1$.

This is precisely the $k=1$ instance of the ``exchange of fixed points across
the Weyl-chamber wall $\lambda_1=\lambda_2$'' generalized in Part~III, and is
the minimal, non-trivial example from which the whole theory of Part~III
takes its start.

\subsection{The Legendre-transformed surface (front) in three dimensions}

So far we have treated the set of degenerate points ($\det H=0$) as a curve
in the $\xi$-plane. Here we change viewpoint and visualize the Legendre
transform itself directly as a three-dimensional surface, to see how the
degenerate points appear on that surface. Since $f$ is non-convex (the
origin is a local maximum), the ordinary single-valued Legendre transform of
convex analysis does not apply. Instead, for each $\xi$ we compute the
gradient $p(\xi)=\nabla f(\xi)$ and the corresponding value
$S^*(\xi):=\xi\cdot p(\xi)-f(\xi)$, and consider the surface (the
``front,'' in the language of singularity theory) swept out as $\xi$ varies:
\[
\Sigma_f=\{(p_1(\xi),p_2(\xi),S^*(\xi)):\xi\in\R^2\}\subset\R^3.
\]
Figure~6 shows the isotropic case. The surface in panel (b) has a
``petal''-like structure near the origin, where several sheets fold over one
another --- a visible manifestation of the non-invertibility (multivaluedness
of the Legendre transform) of $\Lambda_f$.

\begin{figure}[htbp]
\centering
\includegraphics[width=\textwidth]{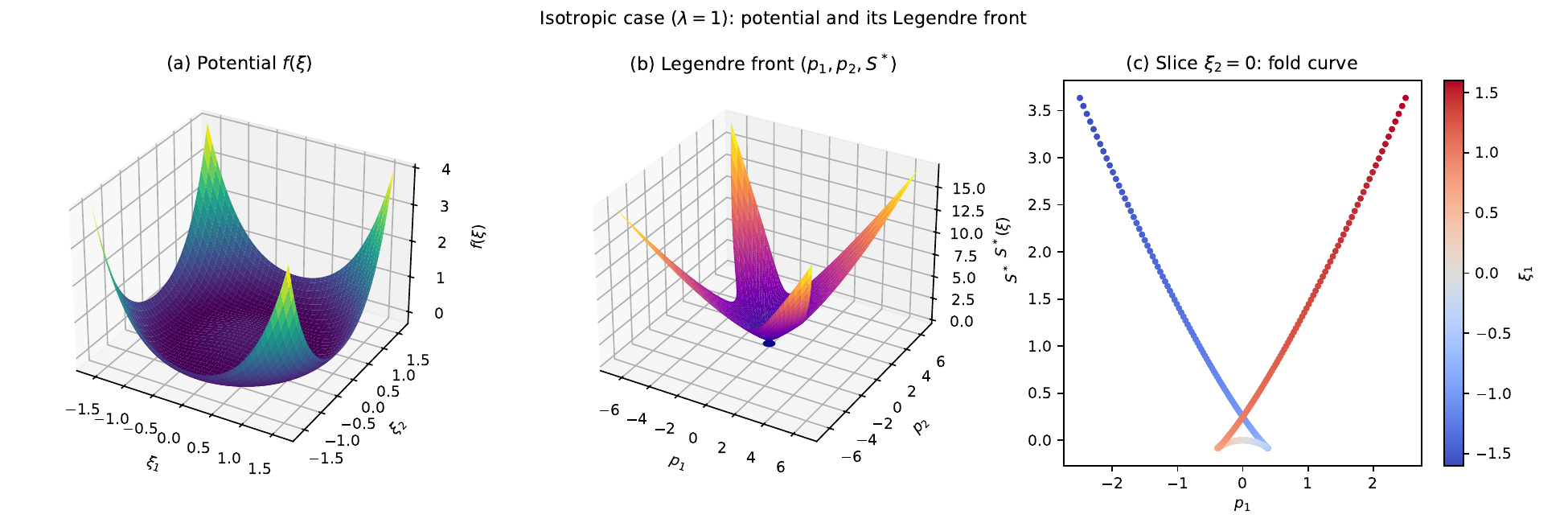}
\caption{Isotropic case ($\lambda=1$): (a) the graph of the potential
$f(\xi)$; (b) the surface $\Sigma_f=(p_1(\xi),p_2(\xi),S^*(\xi))$ after the
Legendre transform; (c) the curve $(p_1,S^*)$ on the slice $\xi_2=0$ (color
indicates the original value of $\xi_1$).}
\end{figure}

Looking more closely at the slice $\xi_2=0$ (Figure~7), the gradient map
$p_1(\xi_1,0)$ is seen to be non-monotonic, forming a peak and a trough at
the caustic locations ($\xi_1=\pm1/\sqrt3,\pm1$). This non-monotonicity is
precisely the origin of the multivaluedness of the Legendre map --- several
$\xi_1$ give the same slope $p_1$ --- and the Maslov index is nothing but a
signed count of how many times a path passes through this fold as it
traverses the curve.

\begin{figure}[htbp]
\centering
\includegraphics[width=\textwidth]{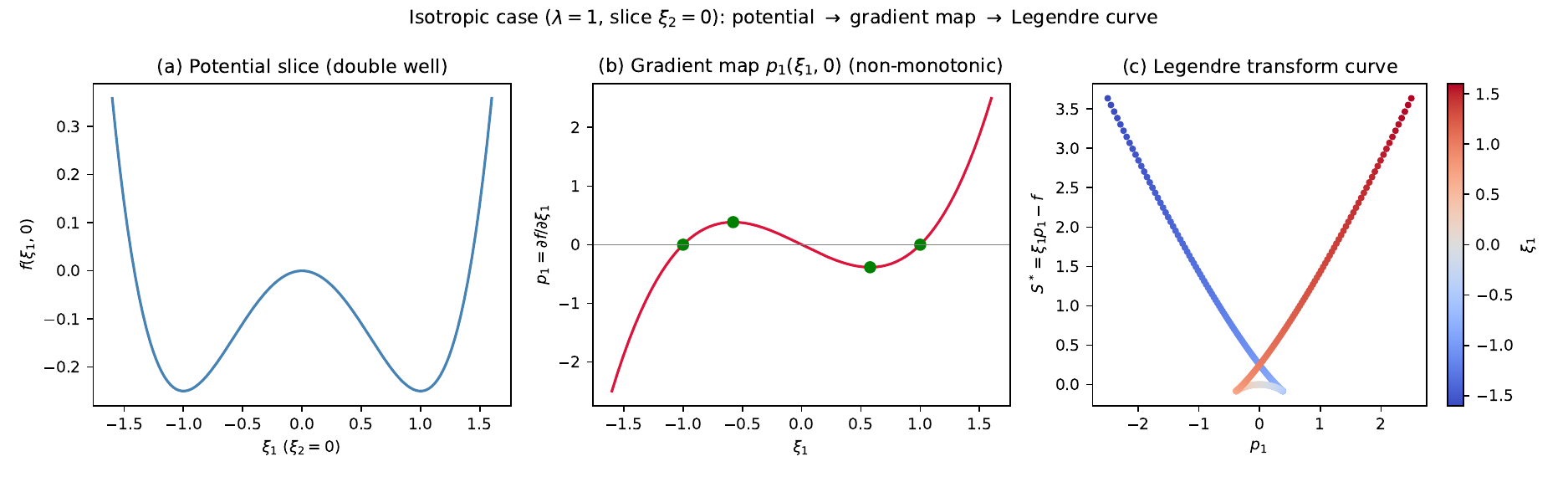}
\caption{Isotropic case ($\lambda=1$, slice $\xi_2=0$): (a) the slice
$f(\xi_1,0)$ of the potential (a double well); (b) the slice $p_1(\xi_1,0)$
of the gradient map (non-monotonic; green points mark the caustic
locations); (c) the corresponding Legendre-transform curve.}
\end{figure}

In the anisotropic case, under the cusp condition $\lambda_2=3\lambda_1$, the
front surface, mapped into $p$-space, becomes a star-shaped (astroid-like)
curve with \textbf{four sharp vertices (cusp points)} (Figure~8). This is a
concrete instance of the phenomenon that Ekeland~\cite{ekeland1977} rigorously
formalized in the general theory of non-convex Legendre transforms: ``a
simple inflection point gives birth to a simple cusp.''

\begin{figure}[htbp]
\centering
\includegraphics[width=\textwidth]{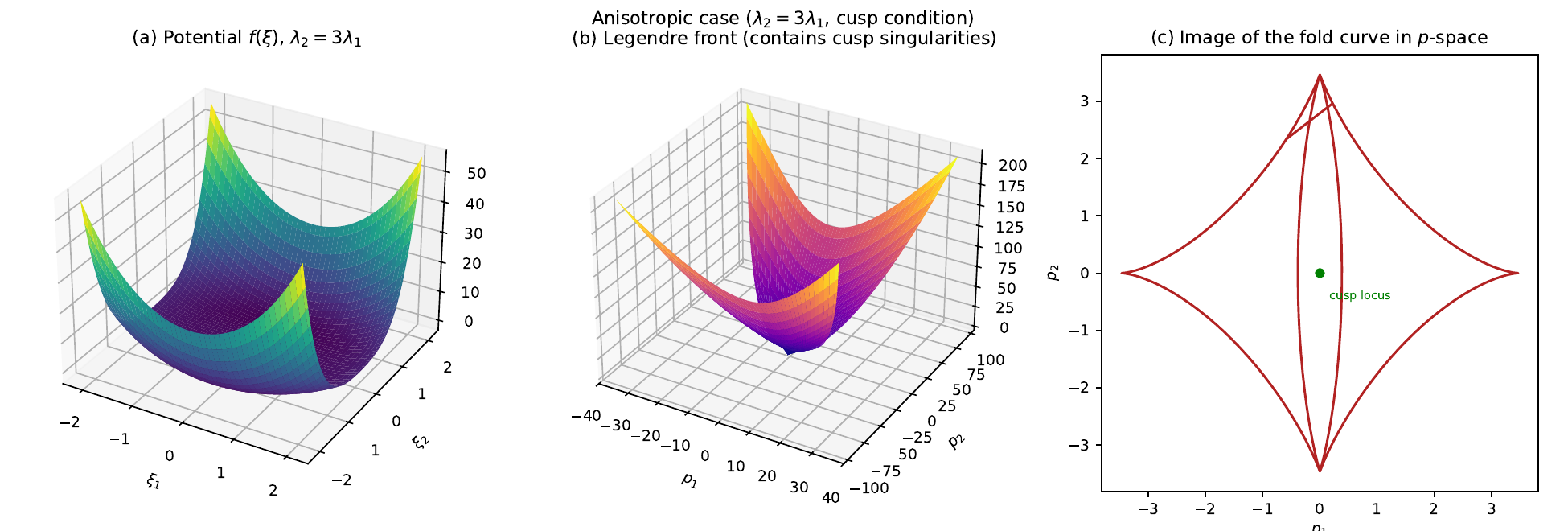}
\caption{Anisotropic case ($\lambda_2=3\lambda_1$, the cusp condition): (a)
the graph of the potential $f$; (b) the surface after the Legendre transform
(containing cusp singularities); (c) the image of the fold curve in
$p$-space (the sharp points are the cusp singularities).}
\end{figure}

\subsection{Explicit computation of the Maslov index}

Along the radial path $\rho:0\to2$, the sign of
$\det H\propto(3\rho^2-1)(\rho^2-1)$ flips twice:
\[
\rho\in\Big(0,\tfrac{1}{\sqrt3}\Big):+1,\qquad
\rho\in\Big(\tfrac{1}{\sqrt3},1\Big):-1,\qquad
\rho\in(1,2):+1,
\]
(Figure~9); the local indices are $\ind_{\rho=1/\sqrt3}=-1$,
$\ind_{\rho=1}=+1$, so the net Maslov index is
$\mathrm{MaslovInd}(0\to2)=0$. This reflects the fact that the start and end
points lie in the same ``chamber'' (sign $+1$); the unsigned number of
crossings is $2$. Along the partial path $\rho:0\to0.8$ (crossing only
$\rho=1/\sqrt3$), $\mathrm{MaslovInd}=-1$, non-zero.

\begin{figure}[htbp]
\centering
\includegraphics[width=0.6\textwidth]{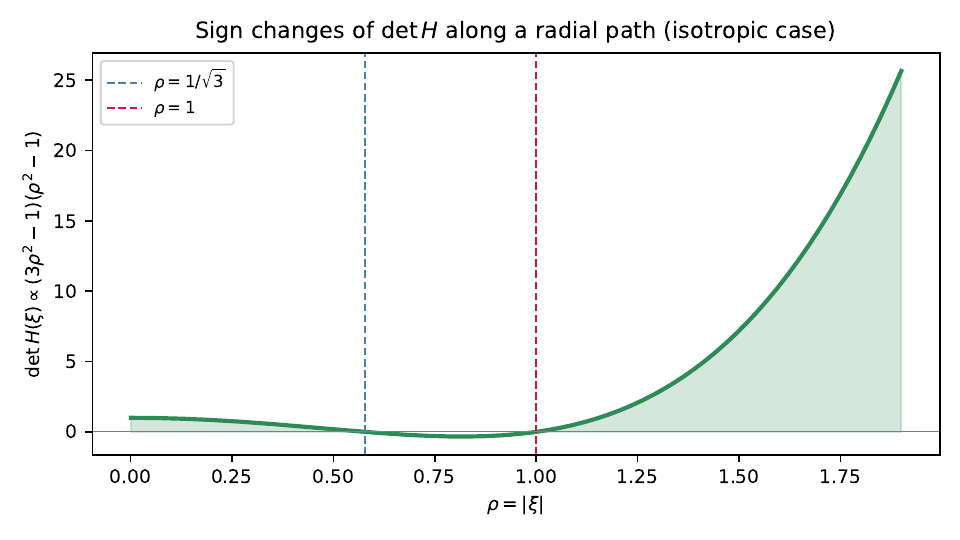}
\caption{Sign changes of $\det H$ and the local Maslov contribution along a
radial path (isotropic case).}
\end{figure}

\begin{remark}[Relation to other deformation families; scope of this paper]
$f$ can be situated within a broader one-parameter deformation family: the
information-geometric Box--Cox-type potential $h_\alpha$, and the pure-power
unified family $\Phi_\alpha$ (with $\Phi_2\equiv f$). We examine both in
detail in Sections~4--6. What is actually needed for Part~III is only
$f=\Phi_2$ (the Oja--Brockett potential with $B=I$), but the comparison of
$h_\alpha$ and $\Phi_\alpha$ has independent value in illuminating a general
principle of the Maslov index --- the distinction between the net index and
the crossing count --- and we survey it below.
\end{remark}

\section{Review of the Series $h_\alpha$}

Consider the potential corresponding to the Box--Cox/Tsallis-type deformed
exponential family of information geometry,
\begin{equation}
h_\alpha(\xi)=\frac12\xi^\top A\xi-\frac{1}{2\alpha}\big[(\xi^\top\xi+1)^\alpha-1\big]. \label{eq:halpha}
\end{equation}
Whereas $f$ has the structure of raising $Q=\xi^\top A\xi$ itself to a power,
$h_\alpha$ introduces an absolute reference length ($+1$) into the Euclidean
norm $\xi^\top\xi$ and deforms its logarithm. In the isotropic case
($A=\lambda I$), the caustic radius is given by the closed formula, depending
on both $\lambda$ and $\alpha$,
\[
\rho_T^2=\lambda^{1/(\alpha-1)}-1,\qquad
\lambda=u^{\alpha-2}\big[(2\alpha-1)u-2(\alpha-1)\big],\quad u=\rho_R^2+1,
\]
and a phase-transition-like reversal of the existence condition occurs ---
degenerate points exist only for $\lambda\ge1$ when $\alpha>1$, and only for
$\lambda\le1$ when $\alpha<1$ --- a phenomenon absent for $f$ (Figures~10, 11).

\textbf{Reversal of the existence condition ($\alpha<1$):} Fixing
$\alpha=0.5$ and comparing $\lambda=2$ (always positive, no degeneracy) with
$\lambda=0.5$ (a sign reversal partway through, degeneracy present)
(Figure~12), we see that the sign of the tangential Hessian component can
become negative only for $\lambda<1$. Tracking the $\lambda$-dependence
(Figure~13), no degenerate point exists for $\lambda>1$, and one is born and
annihilated exactly at the origin when $\lambda=1$. Figure~14 shows how the
caustic circle changes as $\lambda$ is moved below $1$.

\begin{figure}[htbp]
\centering
\includegraphics[width=\textwidth]{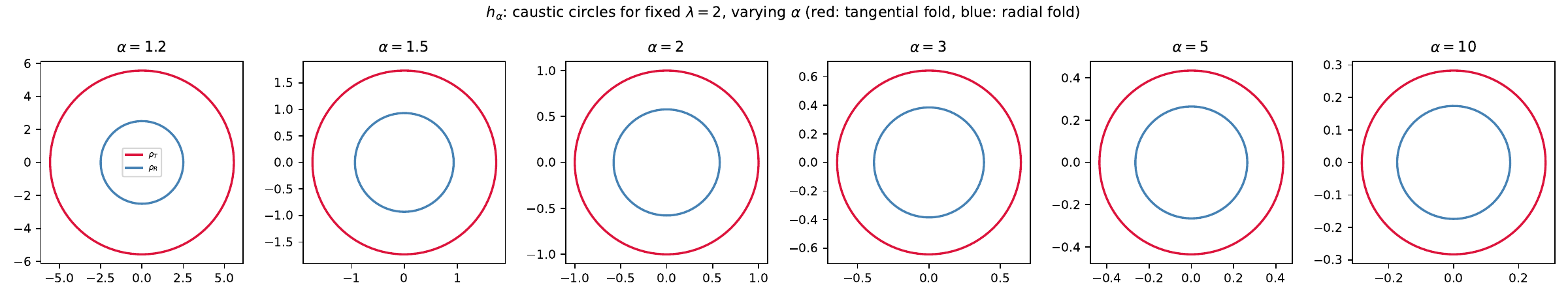}
\caption{The two caustic circles (red: tangential fold, blue: radial fold)
for fixed $\lambda=2$ as $\alpha$ varies.}
\end{figure}

\begin{figure}[htbp]
\centering
\includegraphics[width=0.6\textwidth]{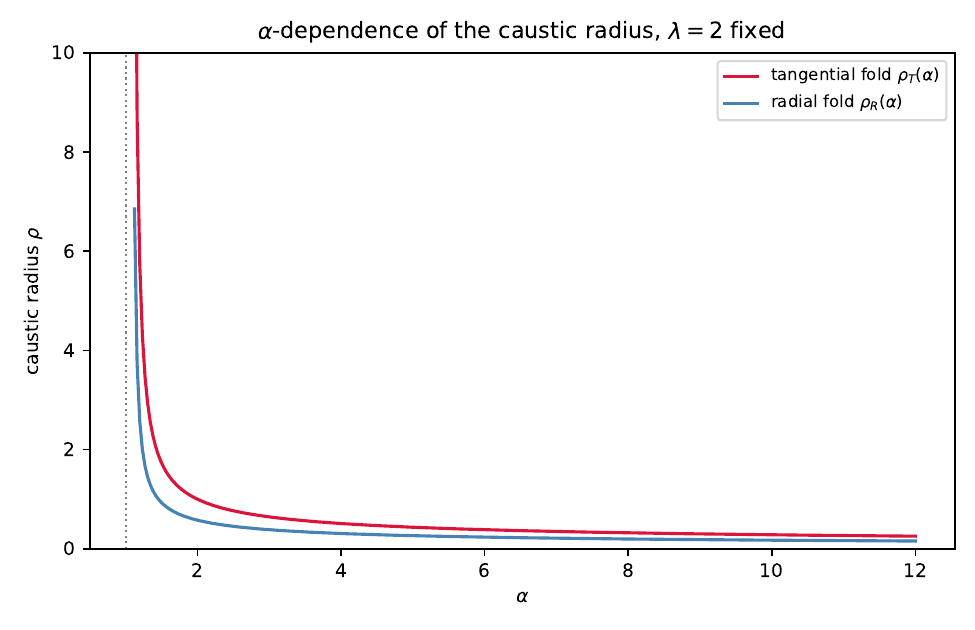}
\caption{$\alpha$-dependence of the caustic radius for fixed $\lambda=2$. As
$\alpha\to1^+$ the circle diverges to infinity.}
\end{figure}

\begin{figure}[htbp]
\centering
\includegraphics[width=\textwidth]{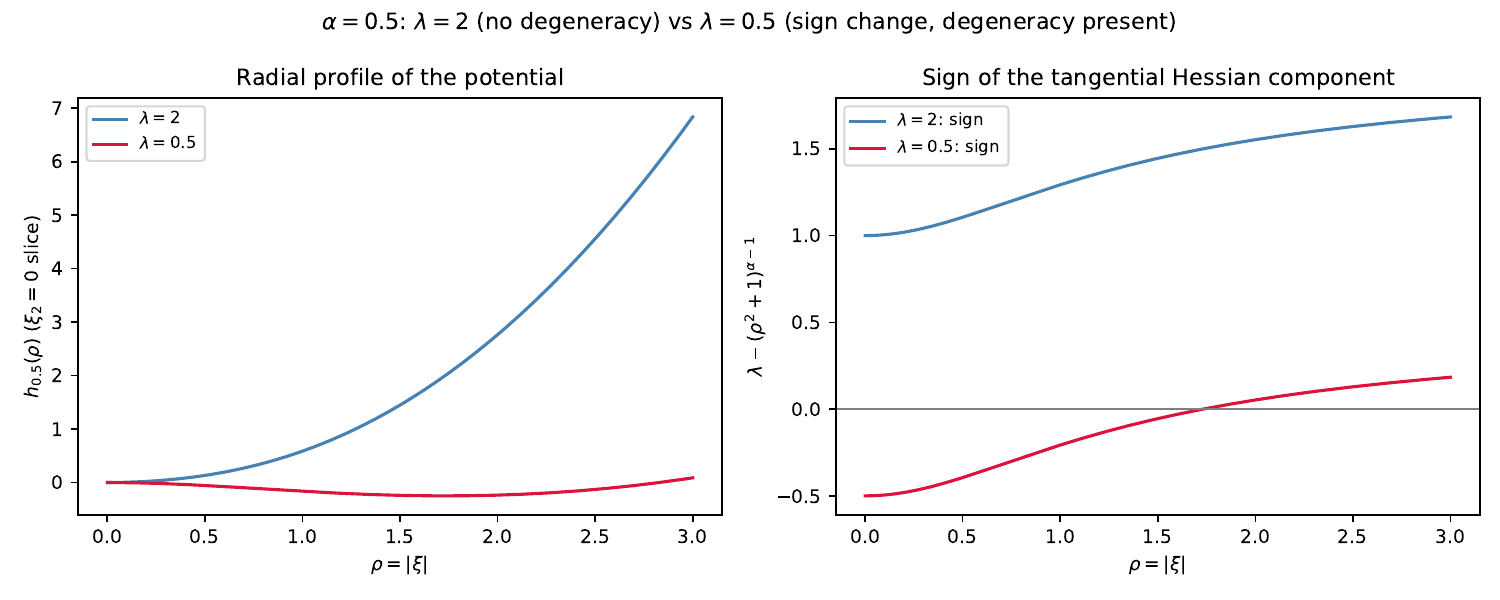}
\caption{$\alpha=0.5$: comparison of $\lambda=2$ (always positive, no
degeneracy) with $\lambda=0.5$ (sign reversal partway through, degeneracy
present). Left: radial slice of the potential. Right: sign of the tangential
Hessian component.}
\end{figure}

\begin{figure}[htbp]
\centering
\includegraphics[width=0.6\textwidth]{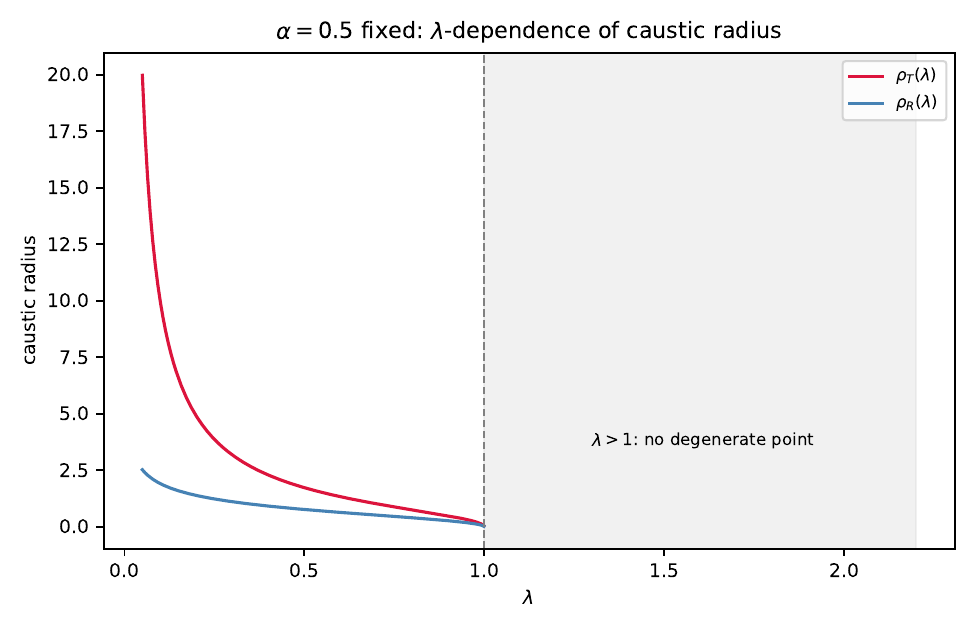}
\caption{$\lambda$-dependence of the caustic radius for fixed $\alpha=0.5$.
A degenerate point is born/annihilated at the origin at $\lambda=1$; none
exists for $\lambda>1$.}
\end{figure}

\begin{figure}[htbp]
\centering
\includegraphics[width=\textwidth]{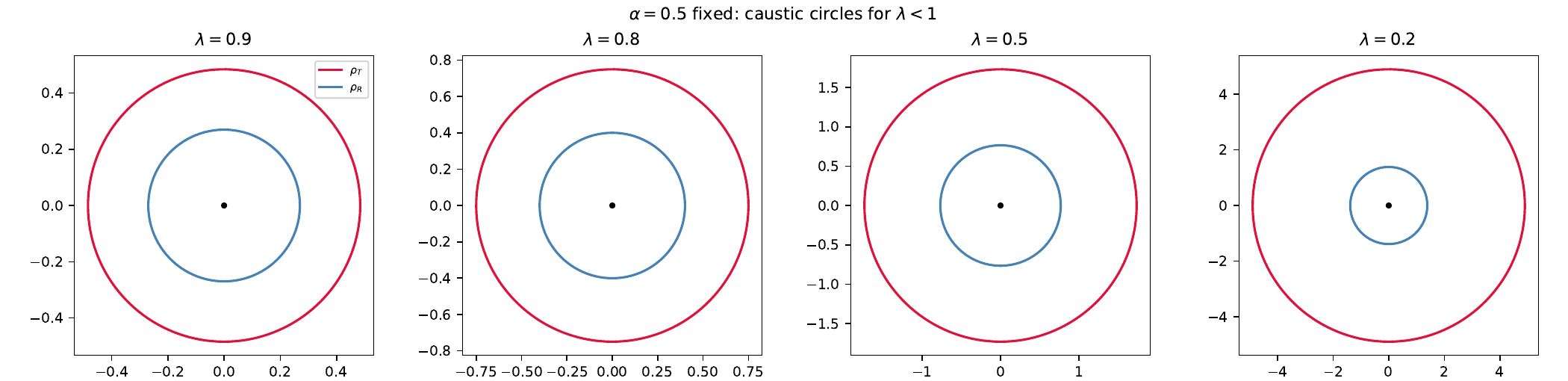}
\caption{Caustic circles for fixed $\alpha=0.5$ as $\lambda$ varies below $1$.}
\end{figure}

\textbf{Path-dependence of the Maslov index:} along a path through the fixed
observation point $\rho_0=0.9$ moving $\lambda$ from $2.0\to0.2$
(Figure~15), the sign flips at the two points $\lambda=0.743,0.411$, and the
net Maslov index is again $0$ (the same structure as for $f$). Figure~16
shows the sign change along the radial path $\rho:0\to3$ (with
$\lambda=0.5,\alpha=0.5$ fixed).

\begin{figure}[htbp]
\centering
\includegraphics[width=0.6\textwidth]{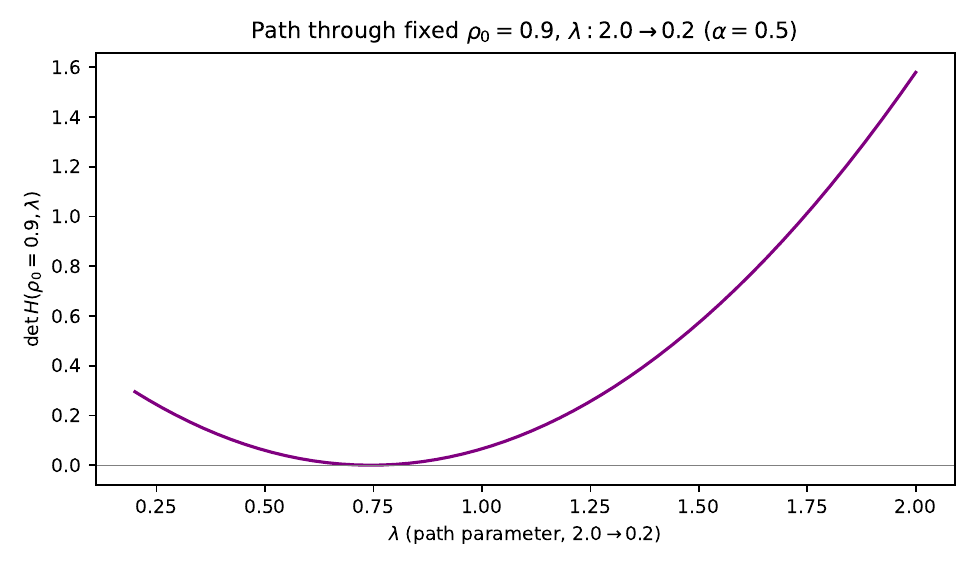}
\caption{Change of $\det H$ along a path through the fixed observation point
$\rho_0=0.9$, moving $\lambda$ from $2.0\to0.2$ ($\alpha=0.5$).}
\end{figure}

\begin{figure}[htbp]
\centering
\includegraphics[width=0.6\textwidth]{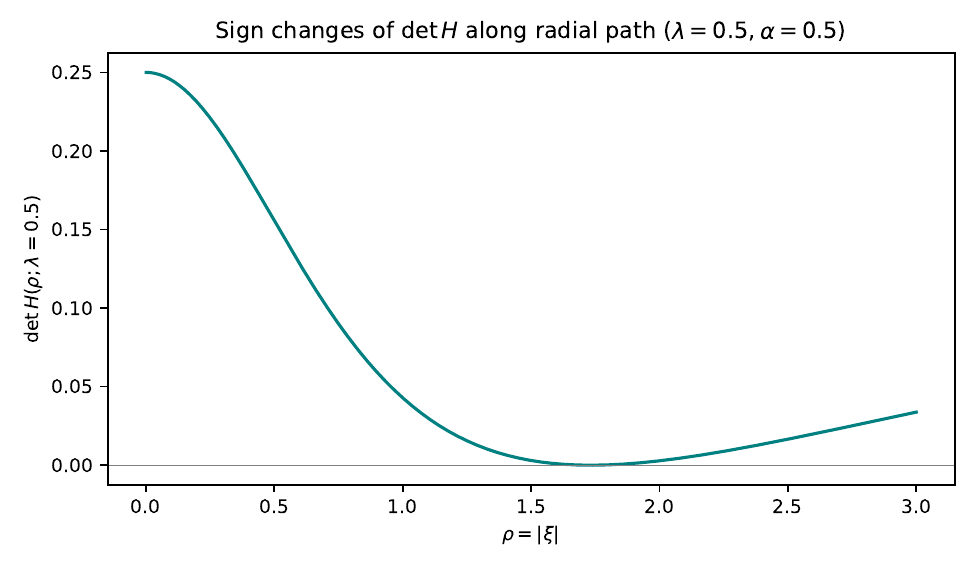}
\caption{Sign change of $\det H$ along the radial path $\rho:0\to3$
($\lambda=0.5,\alpha=0.5$).}
\end{figure}

In the anisotropic case, the cusp condition becomes, at $\alpha=2$, the
affine relation with an intercept, $\lambda_2=\lambda_1/3+2/3$
(Figure~17), a different form from the proportionality $\lambda_2=3\lambda_1$
found for $f$. Figure~18 compares the cusp curves across $\alpha$;
Figure~19 compares them with the case $\alpha=0.5$; Figure~20 shows the
anisotropic deformation at $\alpha=0.5$. Figure~21 compares the sign change
along the radial path for $\alpha=1.5,2,5$ with $\lambda=2$ fixed.

\begin{figure}[htbp]
\centering
\includegraphics[width=\textwidth]{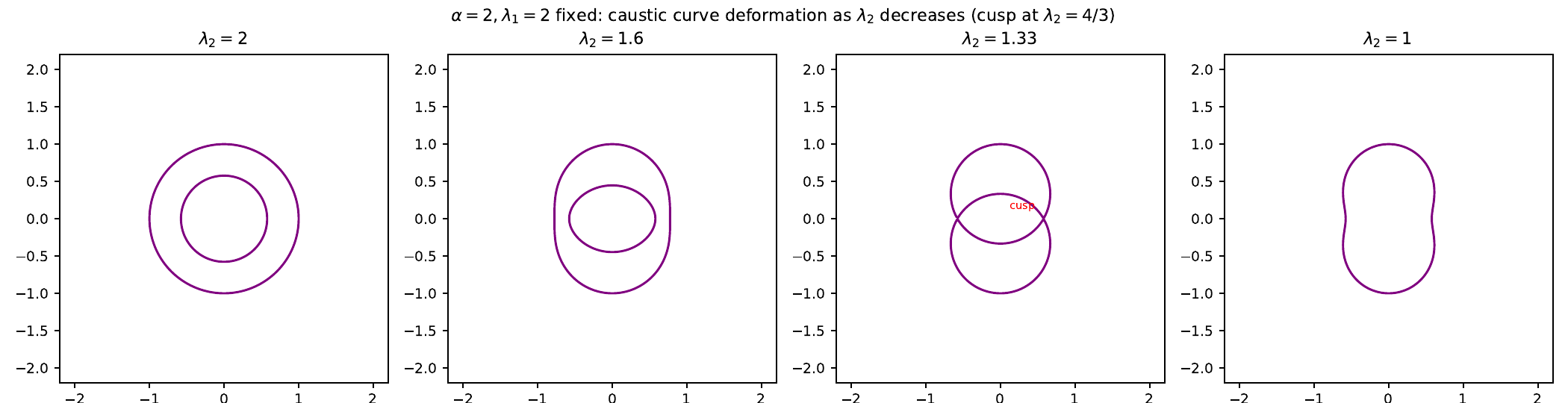}
\caption{Deformation of the caustic curve at $\alpha=2,\lambda_1=2$ fixed,
as $\lambda_2$ decreases; a cusp forms at $\lambda_2=4/3$.}
\end{figure}

\begin{figure}[htbp]
\centering
\includegraphics[width=0.5\textwidth]{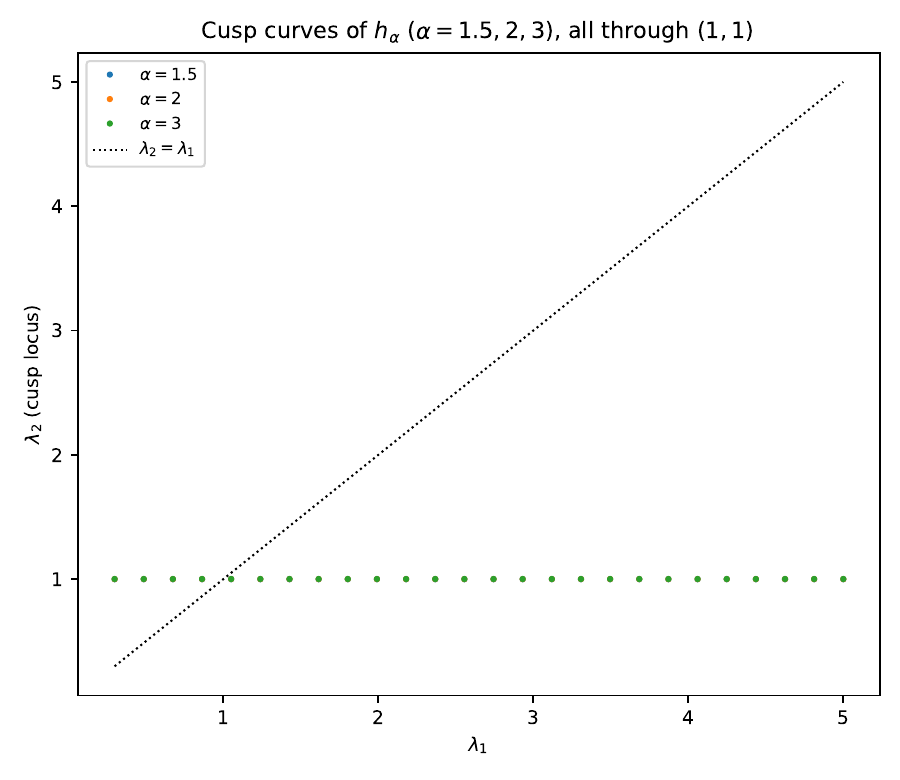}
\caption{Cusp curves of $h_\alpha$ ($\alpha=1.5,2,3$). All pass through $(1,1)$.}
\end{figure}

\begin{figure}[htbp]
\centering
\includegraphics[width=0.5\textwidth]{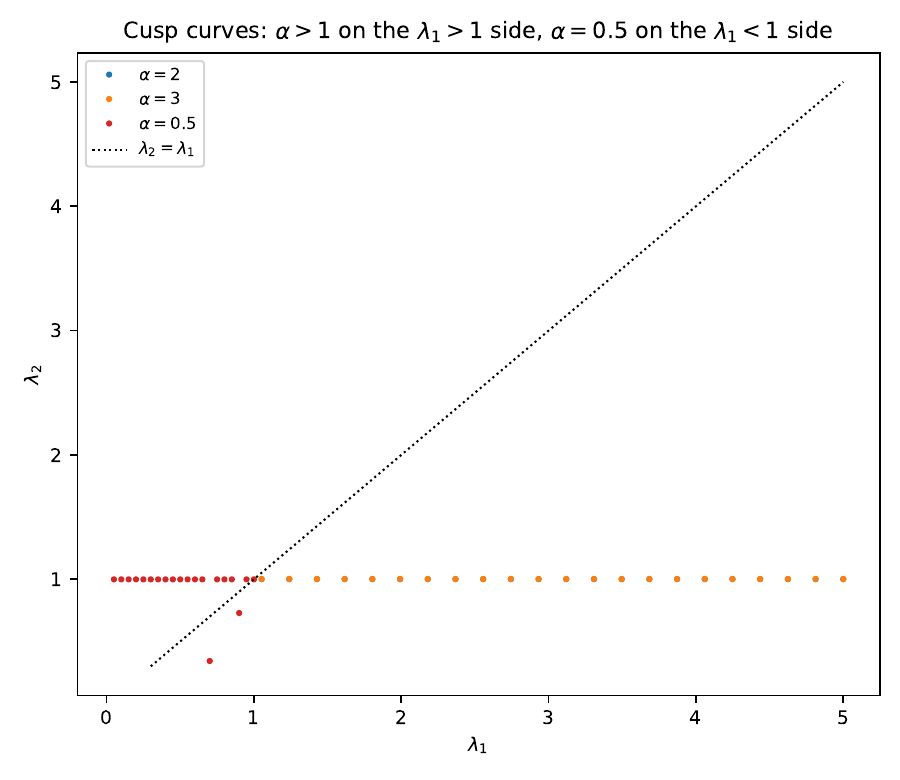}
\caption{Comparison of cusp curves: $\alpha>1$ lies on the $\lambda_1>1$ side,
$\alpha=0.5$ on the $\lambda_1<1$ side.}
\end{figure}

\begin{figure}[htbp]
\centering
\includegraphics[width=\textwidth]{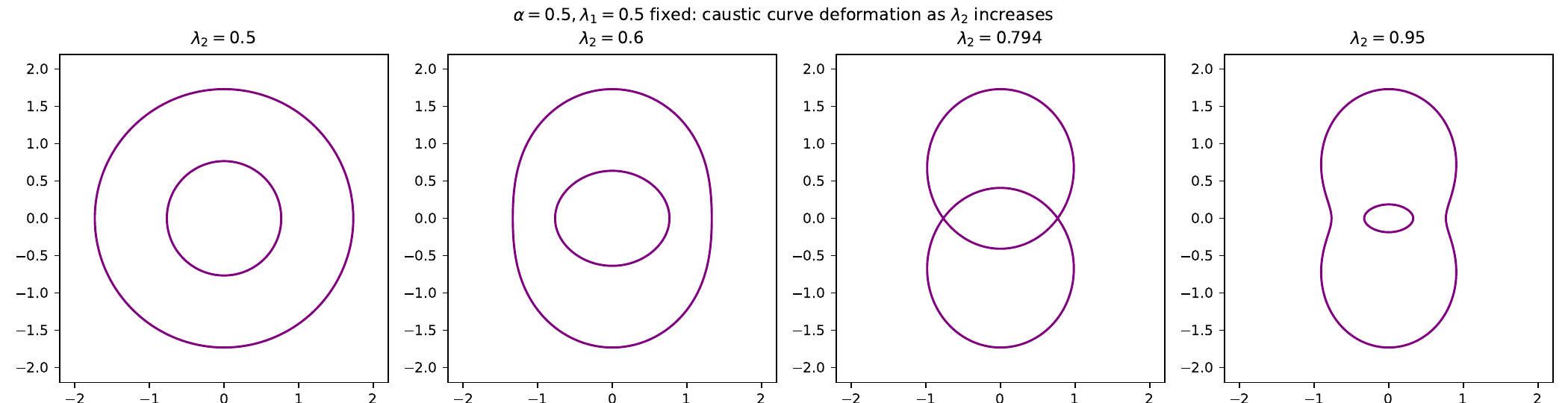}
\caption{Deformation of the caustic curve at $\alpha=0.5,\lambda_1=0.5$
fixed, as $\lambda_2$ increases.}
\end{figure}

\begin{figure}[htbp]
\centering
\includegraphics[width=\textwidth]{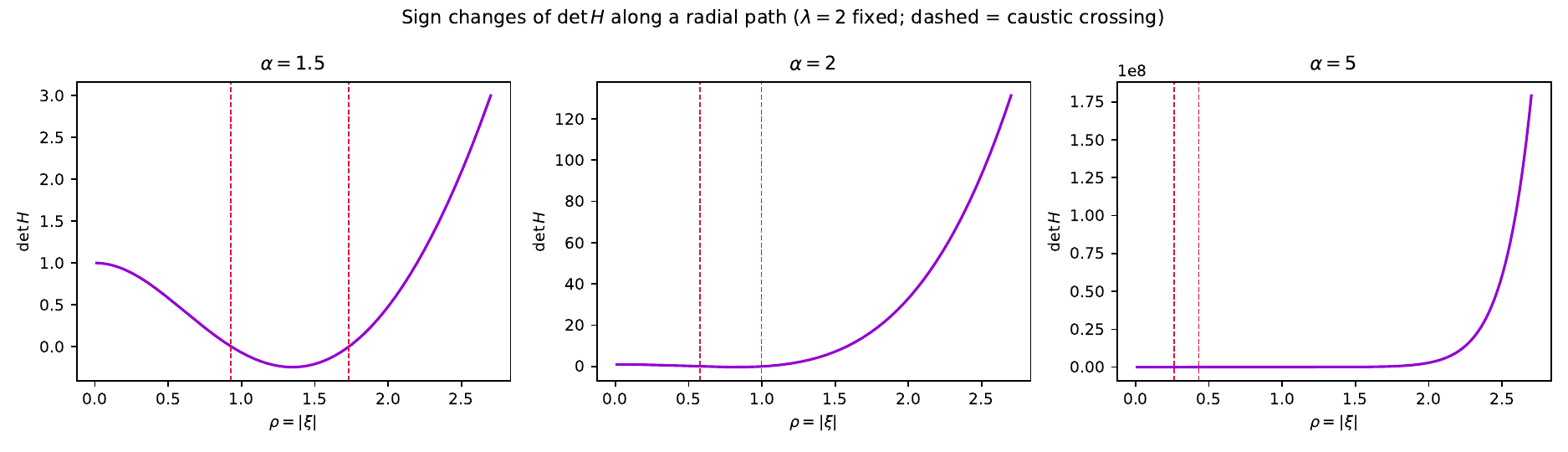}
\caption{Sign change of $\det H$ along a radial path ($\lambda=2$ fixed;
dashed lines mark sign reversals, i.e., caustic crossings).}
\end{figure}

\section{Consistency Check: Does $h_2$ Coincide with $f$?}

Substituting $\alpha=2$ into the formula \eqref{eq:halpha} gives
$h_2(\xi)=\tfrac12Q-\tfrac14[(\xi^\top\xi+1)^2-1]$, which, like $f$, has the
structure of a ``quadratic term plus a quartic term.'' Moreover, at
$\lambda=2$ the two caustic radii agree exactly, $\rho_T=1,\rho_R=1/\sqrt3$.
Is this a coincidence?

\begin{proposition}[Proof of non-coincidence]
$h_2(\xi)$ and $f(\xi)$ do not agree identically for a general
positive-definite symmetric matrix $A$.
\end{proposition}

This can be verified by direct expansion in a computer-algebra system
(SymPy): the quartic term of $h_2$ does not depend on $A$ at all
($-\tfrac14\xi^\top\xi\cdot\xi^\top\xi$, i.e.\ $-r^4/4$), whereas the quartic
term of $f$ carries coefficients quadratic in the entries of $A$.

In the isotropic case, $f$ is written entirely as a multiple of $\lambda^2$,
so the location of its caustic does not depend on $\lambda$ at all; but for
$h_2$ the coefficient of the quadratic term is linear in $\lambda$ while the
quartic term is a $\lambda$-independent constant, so their relative weight
changes with $\lambda$. It is only at the particular value $\lambda=2$ that
these two different dependencies happen to give the same numerical value (an
accidental coincidence). Figure~22 overlays the caustic circles of $h_2$ and
$\Phi_2=f$ (see the next section) for several values of $\lambda$, confirming
that they coincide only at $\lambda=2$.

\begin{figure}[htbp]
\centering
\includegraphics[width=\textwidth]{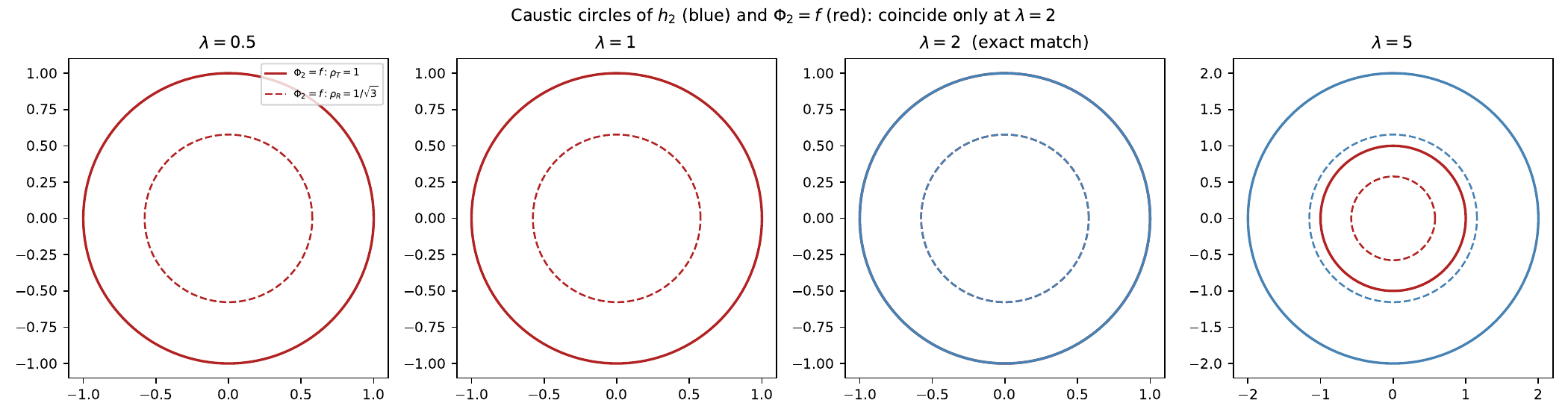}
\caption{Comparison of the caustic circles of $h_2$ (blue) and $\Phi_2=f$
(red, see the next section). At $\lambda=0.5$, $h_2$ has no degenerate point;
at $\lambda=1,5$ the two radii disagree. They coincide only at $\lambda=2$.}
\end{figure}

This check establishes that the conjecture ``$\alpha=2$ of $h_\alpha$ is a
generalization of $f$'' is not correct. So how should a unified family that
correctly contains $f$ be constructed?

\section{Construction of a New Unified Family $\Phi_\alpha$}

The quartic term $\tfrac14Q^2$ of $f$ is written as $Q=\xi^\top A\xi$ raised
to a power. Building on this observation, we introduce a deformation family
of ``pure-power type,'' distinct from the ``Box--Cox type'' of $h_\alpha$:
\begin{equation}
\Phi_\alpha(\xi):=-\frac12\xi^\top A^2\xi+\frac{1}{2\alpha}(\xi^\top A\xi)^\alpha. \label{eq:Phialpha}
\end{equation}

\begin{theorem}[Main theorem of the unified family]
$\Phi_\alpha$ satisfies: (i) $\Phi_2(\xi)\equiv f(\xi)$ identically, for any
positive-definite symmetric $A$; (ii) $\Phi_1(\xi)=\tfrac12\xi^\top
A(I-A)\xi$ degenerates to a pure quadratic form; (iii) in the isotropic case
the caustic radii are $\rho_T^2=\lambda^{(2-\alpha)/(\alpha-1)}$,
$\rho_R^2=\rho_T^2/(2\alpha-1)^{1/(\alpha-1)}$; (iv) in the anisotropic case
the cusp condition is $\lambda_2=(2\alpha-1)\lambda_1$ (a proportionality
through the origin).
\end{theorem}

Figure~23 shows the isotropic potential shape of $\Phi_\alpha$; $\alpha=2$
coincides exactly with $f$ (a Mexican hat).

\begin{figure}[htbp]
\centering
\includegraphics[width=\textwidth]{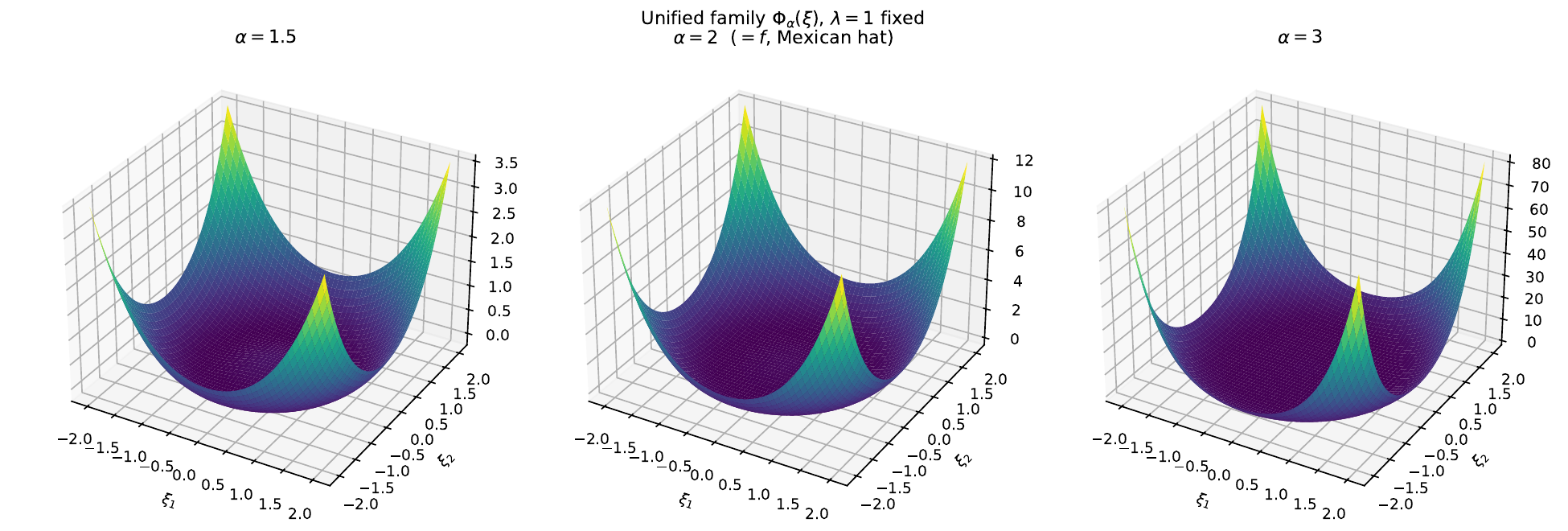}
\caption{Shape of the unified family $\Phi_\alpha(\xi)$ for $\alpha=1.5,2,3$
($\lambda=1$ fixed).}
\end{figure}

The exponent $(2-\alpha)/(\alpha-1)$ in (iii) vanishes exactly at $\alpha=2$;
this is the first value at which both terms of $\Phi_\alpha$ scale by the
same power $c^2$ under $A\to cA$ (the whole function becomes homogeneous).
This is the true reason for the special scale invariance of $f=\Phi_2$
(Figure~24).

\begin{figure}[htbp]
\centering
\includegraphics[width=0.6\textwidth]{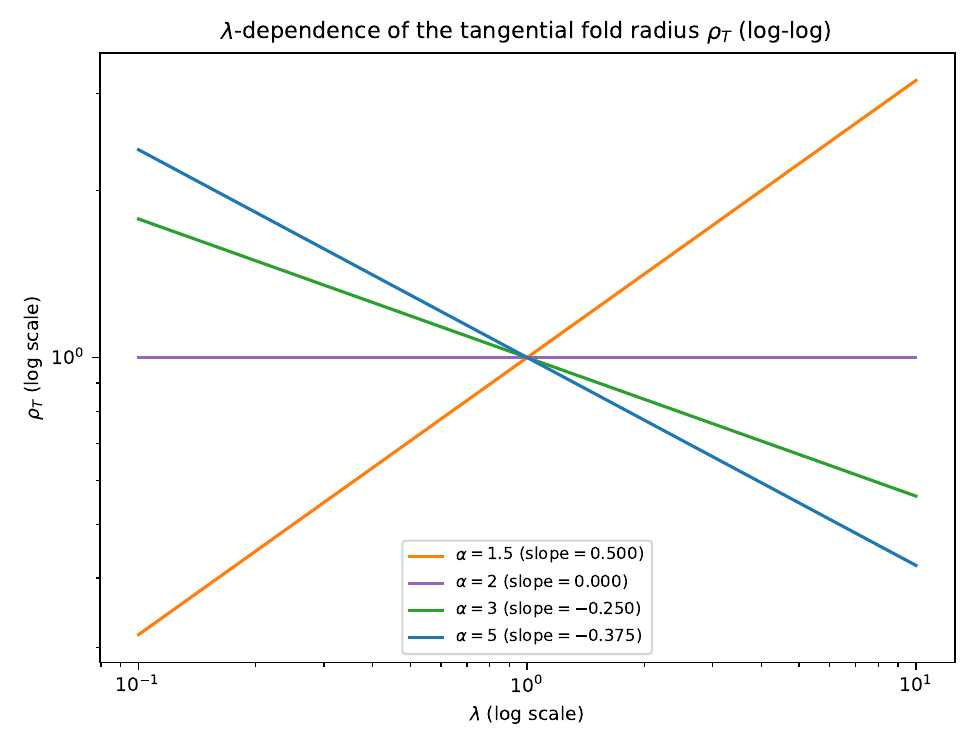}
\caption{$\lambda$-dependence of the tangential fold radius $\rho_T$ of
$\Phi_\alpha$ (log-log plot). The slope $(2-\alpha)/(2(\alpha-1))$ is exactly
$0$ at $\alpha=2$ (horizontal line), recovering the scale invariance of $f$.}
\end{figure}

The cusp condition $\lambda_2=(2\alpha-1)\lambda_1$ of (iv) becomes, exactly
at $\alpha=2$, $\lambda_2=3\lambda_1$, in precise agreement with the result
for $f$. In contrast to the cusp condition of $h_\alpha$ (which carries an
intercept), this remains a proportionality through the origin for every
value of $\alpha$.

\section{Structural Comparison of the Two Deformation Mechanisms}

Both $h_\alpha$ and $\Phi_\alpha$ share the common spirit of ``deforming a
quadratic form nonlinearly, taking $\alpha=1$ or $\alpha=2$ as a reference
point,'' yet the manner of deformation is fundamentally different: $h_\alpha$
introduces an absolute reference length into the Euclidean norm
(Box--Cox/Tsallis type), while $\Phi_\alpha$ raises the quadratic form itself
to a power (pure-power type). This difference gives rise to entirely
different kinds of special behavior: the phase-transition-like reversal of
the existence condition in $h_\alpha$, and the recovery of scale invariance
at $\alpha=2$ in $\Phi_\alpha$. Figure~25 overlays the cusp curves of both:
one sees clearly the structural difference between $\Phi_\alpha$ (a family
of lines through the origin) and $h_\alpha$ (a family of nonlinear curves
through $(1,1)$).

\begin{figure}[htbp]
\centering
\includegraphics[width=0.68\textwidth]{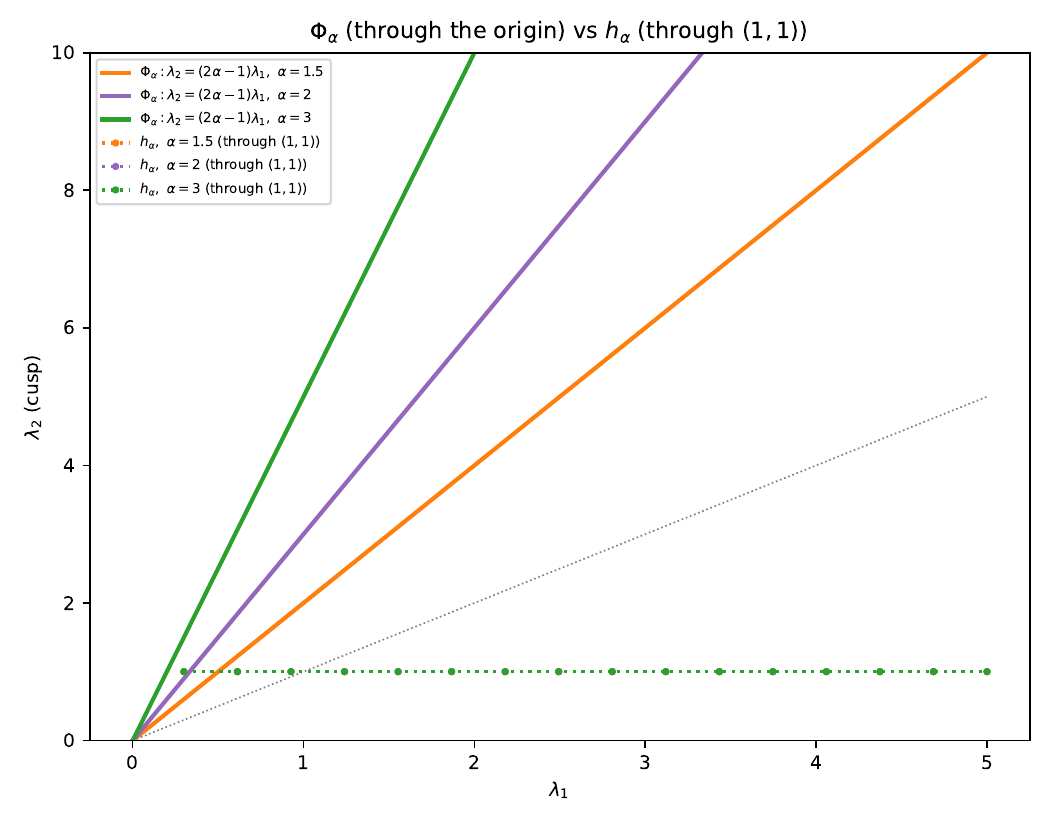}
\caption{Comparison of cusp curves. $\Phi_\alpha$ (solid) is always a family
of lines through the origin, while $h_\alpha$ (dotted) is a family of
nonlinear curves through $(1,1)$.}
\end{figure}

\begin{remark}[Notation used from here on]
From this point on we use only $f=\Phi_2$ (the Oja--Brockett potential with
$B=I$). The detailed mathematical structure of $h_\alpha,\Phi_\alpha$ has
independent value as a way of confirming, from a different angle, the
general principle of the Maslov index (the distinction between the net
index and the crossing count), but the generalization to matrices from
Part~II onward deals exclusively with the potential of $f$-type,
$\tr[(AXBX^\top)^2]$.
\end{remark}

\part{Generalization to Matrix-Valued Potentials: Numerical Observations}
\begin{center}
\textit{Brockett's Sorting Theorem, Flag Manifolds, and Wall-Crossing in the Weyl Group}
\end{center}

In Part~I we studied, for the potential $f$ on the two-dimensional vector
$\xi\in\R^2$, the relationship between degenerate points of the Legendre map
and the Maslov index. In this part we replace $\xi$ by an $n\times k$ tall
real matrix $X$ and dig, through numerical experiments, into the concrete
connections with mathematics spanning classical control theory and Lie-group
representation theory: Brockett's sorting theorem, flag manifolds, and
wall-crossing in the Weyl group. Every phenomenon observed here is proved
rigorously as a theorem in Part~III --- this part is, in effect, an
exploratory preparatory stage in which we numerically map out ``what Part
III's theory needs to explain.''

\section{Introduction: Generalizing from Two-Dimensional Vectors to Matrices}

We replace $\xi\in\R^2$ by an $n\times k$ tall real matrix $X$ ($n\ge k$),
and, in addition to a positive-definite symmetric matrix $A$ ($n\times n$),
introduce a diagonal positive-definite matrix
$B=\diag(b_1,\dots,b_k)$, $b_1>b_2>\cdots>b_k>0$. When $k=1$, $X$ is simply an
$n$-dimensional vector and $B$ is a $1\times1$ scalar, so the discussion of
Part~I ($n=2,k=1,B=1$) is contained completely as the simplest special case
of this generalization. Once $k\ge2$, two new phenomena appear that were
absent before:
\begin{itemize}
\item \textbf{Sorting:} at the fixed points of the negative gradient flow of
an appropriate potential, the columns of $X$ align, in a form that sorts
$k$ of the eigenvectors of $A$ having the largest eigenvalues in the order
determined by the diagonal entries of $B$ (Brockett's sorting theorem, 1991).
\item \textbf{Convergence to a flag manifold:} when $A$ or $B$ has degenerate
eigenvalues, the fixed points spread from isolated points into a flag
manifold, carrying the continuous symmetry of rotation within the degenerate
subspace.
\end{itemize}
Moreover, if we move the eigenvalues of $A$ continuously through the instant
at which two eigenvalues cross, this ``sort order'' exchanges
discontinuously. This crossing point is a higher-dimensional generalization
of the Maslov cycle (degenerate point) discussed in Part~I.

\section{Setup: Generalizing Four Families of Potentials to Matrices}

Let $X\in\R^{n\times k}$, $A$ be $n\times n$ positive-definite symmetric, and
$B=\diag(b_1,\dots,b_k)$ be a $k\times k$ positive-definite diagonal matrix.
We consider the following four potentials:
\begin{align}
f(X)&=-\frac12\tr(A^2XB^2X^\top)+\frac14\tr[(AXBX^\top)^2], \label{eq:fX}\\
f_\alpha(X)&=-\frac12\tr(A^2XB^2X^\top)+\frac{1}{2\alpha}\tr[(AXBX^\top)^\alpha], \label{eq:falphaX}\\
h_\alpha(X)&=-\frac12\tr(AXBX^\top)+\frac{1}{2\alpha}\tr[(X^\top X+B)^\alpha-I], \label{eq:halphaX}\\
\varphi_\alpha(X)&=\frac12\tr(AXBX^\top)+\frac{\alpha}{4}\|B-X^\top X\|_F^2. \label{eq:phialphaX}
\end{align}
At $k=1$, \eqref{eq:fX} reduces to $f(\xi)$ (up to a factor of $b^2$),
\eqref{eq:falphaX} to $\Phi_\alpha$, and \eqref{eq:halphaX} to $h_\alpha$.
The penalty-type potential $\varphi_\alpha$ of \eqref{eq:phialphaX} is newly
introduced in the matrix setting. It is precisely $f(X)$ of \eqref{eq:fX}
that is the general $(A,B)$ Oja--Brockett potential $V_{A,B}(X)$ rigorously
analyzed in Part~III (using the standard notation there), and this bridges
the numerical observations of this part with the theory of Part~III.

\begin{proposition}[Gradient formulas]
Writing $Y:=AXBX^\top$, $W_h:=X^\top X+B$, $W_\varphi:=B-X^\top X$,
\begin{align}
\nabla f(X)&=YAXB-A^2XB^2, \\
\nabla h_\alpha(X)&=-AXB+XW_h^{\alpha-1}, \\
\nabla\varphi_\alpha(X)&=AXB+\alpha X(X^\top X-B).
\end{align}
\end{proposition}
These have been verified numerically against finite-difference gradients,
with maximum absolute error below $10^{-9}$ (checked for $n=4,k=2$ and
$n=5,k=3$).

\section{Brockett's Sorting Theorem}

R.\,W.\ Brockett showed in his 1991 paper~\cite{brockett1991} (a conference
version appeared already in 1988~\cite{brockett1988}) that the gradient flow
(the double-bracket flow) of the function $\tr(NX^\top AX)$ ($N$ diagonal) on
the Stiefel manifold ($X^\top X=I_k$) converges, almost surely from a generic
initial condition, so that the columns of $X$ become eigenvectors of $A$,
sorted so that the order of the diagonal entries of $N$ matches the order of
the corresponding eigenvalues of $A$. The potentials \eqref{eq:fX}--\eqref{eq:phialphaX}
all contain at their core the same structure as Brockett's function,
$\tr(AXBX^\top)=\tr(BX^\top AX)$.

\textbf{Numerical verification:} with $n=5,k=3$, $A=\diag(5,4,3,2,1)$,
$B=\diag(3,2,1)$, numerically integrating $\dot X=-\nabla f(X)$ from a random
initial condition, the diagonal entries of $X^\top AX$ converge, at an
exponential rate, to $(5,4,3)$ (the top three eigenvalues of $A$), while
$X^\top X$ converges to the identity matrix. That is,
\[
X_*^\top X_*=I_k,\qquad X_*^\top AX_*=\diag(\lambda_{(1)},\dots,\lambda_{(k)}),
\]
with $\lambda_{(1)}>\cdots>\lambda_{(k)}$ the top $k$ eigenvalues of $A$,
confirmed to high precision (Figure~26). This ``axis-aligned'' equilibrium
is exactly the standard form $X_*=[e_{p(1)},\dots,e_{p(k)}]$ treated
rigorously in Part~III (\S15.2).

\begin{figure}[htbp]
\centering
\includegraphics[width=\textwidth]{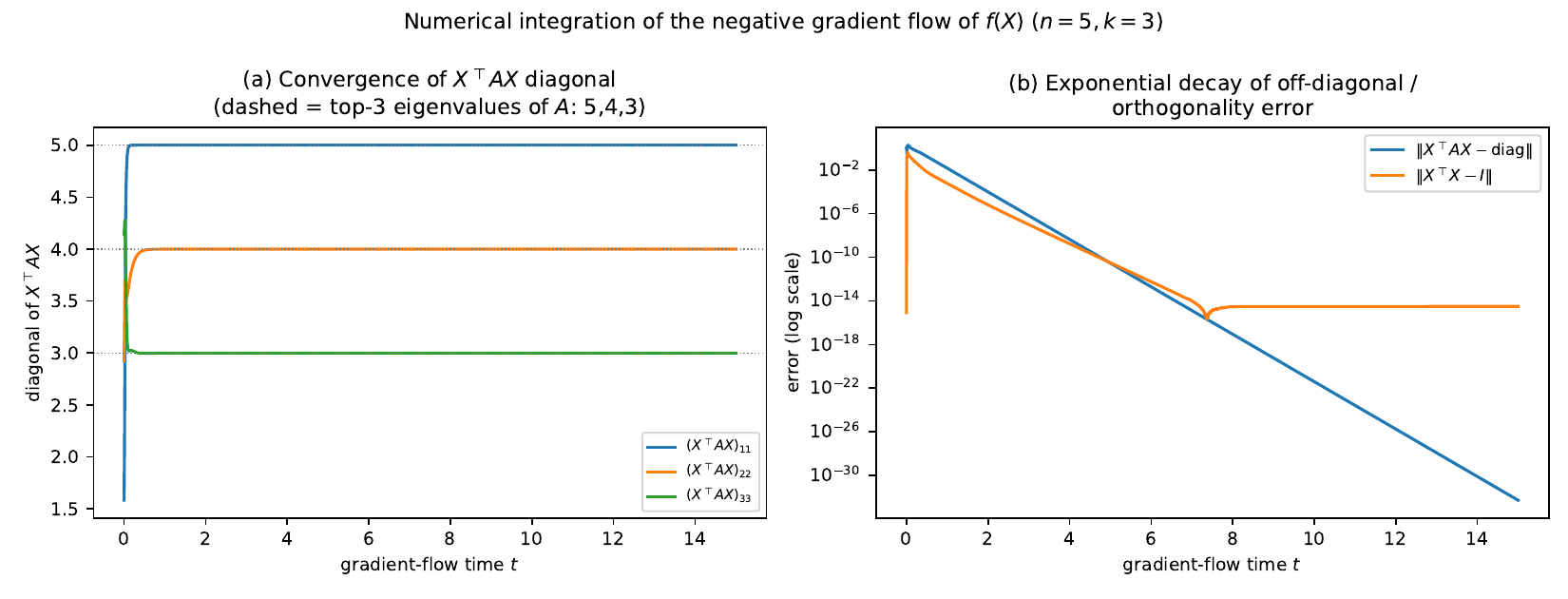}
\caption{Numerical integration of the negative gradient flow of the Type-(1)
potential $f(X)$ ($n=5,k=3$). (a) The diagonal of $X^\top AX$ converges to
the top three eigenvalues of $A$, $(5,4,3)$. (b) The off-diagonal entries and
$X^\top X-I$ decay exponentially.}
\end{figure}

\section{Eigenvalue Degeneracy and Convergence to a Flag Manifold}

When the eigenvalues of $A$ and $B$ are all distinct, the fixed point picked
out by Brockett's sorting theorem is (up to sign flips) isolated. But when
either has a degeneracy (repeated eigenvalue), the situation changes
fundamentally. Since any orthonormal basis of the eigenspace remains a valid
set of eigenvectors, the fixed point is no longer isolated but spreads out
into a quotient space of an orthogonal group (a flag manifold, or
Grassmannian).

\textbf{Degeneracy of $A$'s eigenvalues:} for $n=5,k=3$,
$A=\diag(5,4,4,2,1)$ (eigenvalue $4$ with multiplicity $2$), $B=\diag(3,2,1)$,
solving for the fixed point gives $X_*^\top AX_*=\diag(5,4,4)$,
$X_*^\top X_*=I_3$, and the Hessian on the full space $\R^{15}$ has exactly
one zero eigenvalue (the remaining $14$ are positive). Rotating the two
degenerate columns together through an angle $\theta$, the gradient vanishes
identically for all $\theta$ (to machine precision, $10^{-14}$), and the
potential value is likewise constant.

\textbf{Degeneracy of $B$'s eigenvalues:} symmetrically, keeping
$A=\diag(5,4,3,2,1)$ (distinct) and setting $B=\diag(2,2,1)$ (with
$b_1=b_2$ degenerate), at the fixed point the leading $2\times2$ block
$X_*^\top AX_*$ is \emph{not} diagonalized (off-diagonal entries appear), but
its trace is $4.295+4.705=9=5+4$ (the sum of $A$'s top two eigenvalues), and
the subspace spanned by the first two columns coincides with $A$'s top
two-dimensional eigenspace. Here too the Hessian has exactly one zero
eigenvalue. In both cases the gradient norm and potential value are shown, as
functions of the rotation angle within the degenerate subspace, in
Figure~27: both are identically constant, independent of $\theta$.

\begin{figure}[htbp]
\centering
\includegraphics[width=\textwidth]{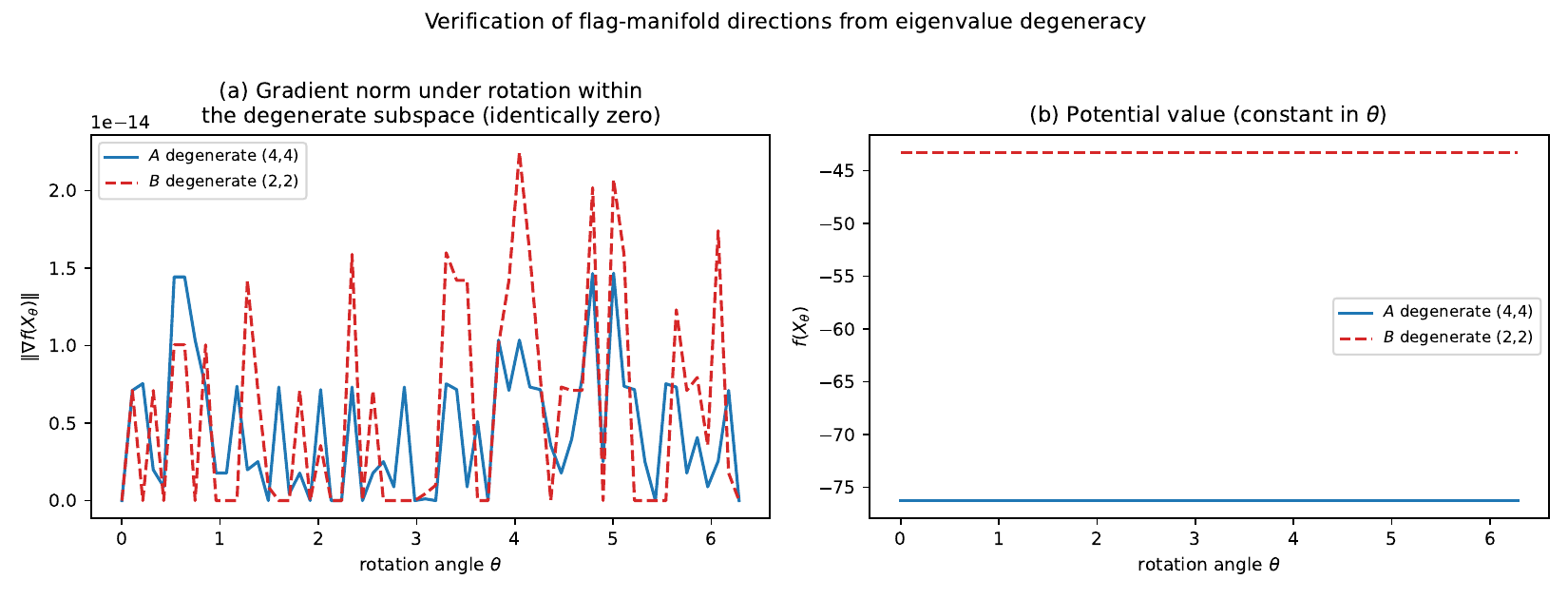}
\caption{Verification of flag-manifold directions from eigenvalue
degeneracy. (a) Gradient norm under rotation of columns within the
degenerate subspace (identically zero in both cases). (b) Corresponding
potential value (constant, independent of the rotation angle).}
\end{figure}

\begin{proposition}[Dimension of the fixed-point set]
Suppose the eigenvalues of $A$ decompose the eigenspace as
$V_1,V_2,\dots$ with multiplicities $(m_1,m_2,\dots)$, and the chosen $k$
columns are taken $k_i$ at a time from each $V_i$. Then each connected
component of the fixed-point set is diffeomorphic to a (partial) flag
manifold of the form $\prod_i\Gr(k_i,m_i)$, of dimension $\sum_ik_i(m_i-k_i)$.
\end{proposition}

\begin{remark}[Connection to Part~III]
The two numerical observations of this section (a one-dimensional zero
direction from degeneracy of $A$, and one from degeneracy of $B$) look at
first like separate phenomena. Lemma~17.1 of Part~III, however, rigorously
shows that both are special cases of a single formula,
\[
\mu^V_{ij}=(\lambda_i-\lambda_j)(b_i-b_j),
\]
(with, respectively, $\lambda_i=\lambda_j$ or $b_i=b_j$). This unifying
structure is not yet visible in this section --- filling that gap is
precisely the theoretical work of Part~III.
\end{remark}

\section{Wall-Crossing in the Weyl Chamber and the Maslov Index}

\subsection{The space of eigenvalue configurations and the Weyl chamber}

The symmetric group $S_n$ acts on the space $\R^n$ with coordinates
$(a_1,\dots,a_n)$ (the eigenvalues of $A$) by permuting coordinates, and the
arrangement of hyperplanes $\{a_i=a_j\}$ where two eigenvalues coincide (the
$A_{n-1}$-type Coxeter arrangement, the \textbf{walls of the Weyl chamber})
divides $\R^n$ into $n!$ congruent chambers. Figure~28 shows the case
$n=3$: three walls divide the space into $6=3!$ chambers (the six possible
orderings of the eigenvalues).

\begin{figure}[htbp]
\centering
\includegraphics[width=0.6\textwidth]{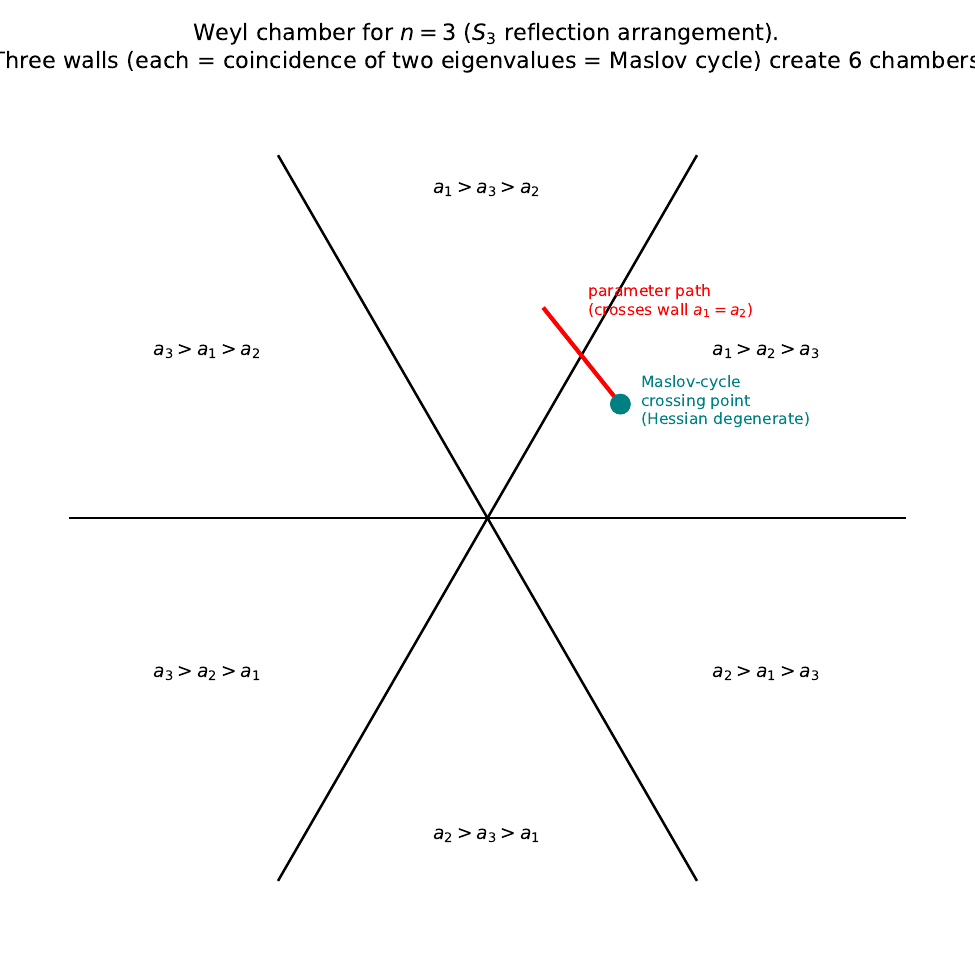}
\caption{The Weyl chamber for $n=3$ (the $S_3$ reflection arrangement).
Three walls (each a coincidence of two eigenvalues = Maslov cycle) divide
the space into six chambers (the six orderings of the eigenvalues). The red
path crosses the wall $a_1=a_2$; the green point marks passage through the
Maslov cycle (Hessian degeneracy).}
\end{figure}

\subsection{Wall-crossing and Hessian degeneracy: numerical verification}

With $n=4,k=2$, $B=\diag(2,1)$, and $A=\diag(a_1,a_2,2,1)$, fixing
$a_1=4$ and sweeping $a_2$ from $3.3$ to $4.7$, the smallest eigenvalue of
the Hessian at the fixed point $X_*(a_2)$ touches zero exactly at
$a_2=4$ (the center of the wall), and returns to positive on either side
(Figure~29(a)). Simultaneously, whereas for $a_2<4$ column~1 tracked
$a_1=4$ (constant) and column~2 tracked $a_2$, for $a_2>4$ the column
correspondence swaps: column~1 always tracks $\max(a_1,a_2)$, column~2
always tracks $\min(a_1,a_2)$ (Figure~29(b)).

This is the exact generalization of the eigenvalue-exchange phenomenon
studied in \S3.2 for the case $n=2,k=1$ --- the bifurcation in which minima
and saddles exchange across $\lambda_1=\lambda_2$.

\begin{figure}[htbp]
\centering
\includegraphics[width=\textwidth]{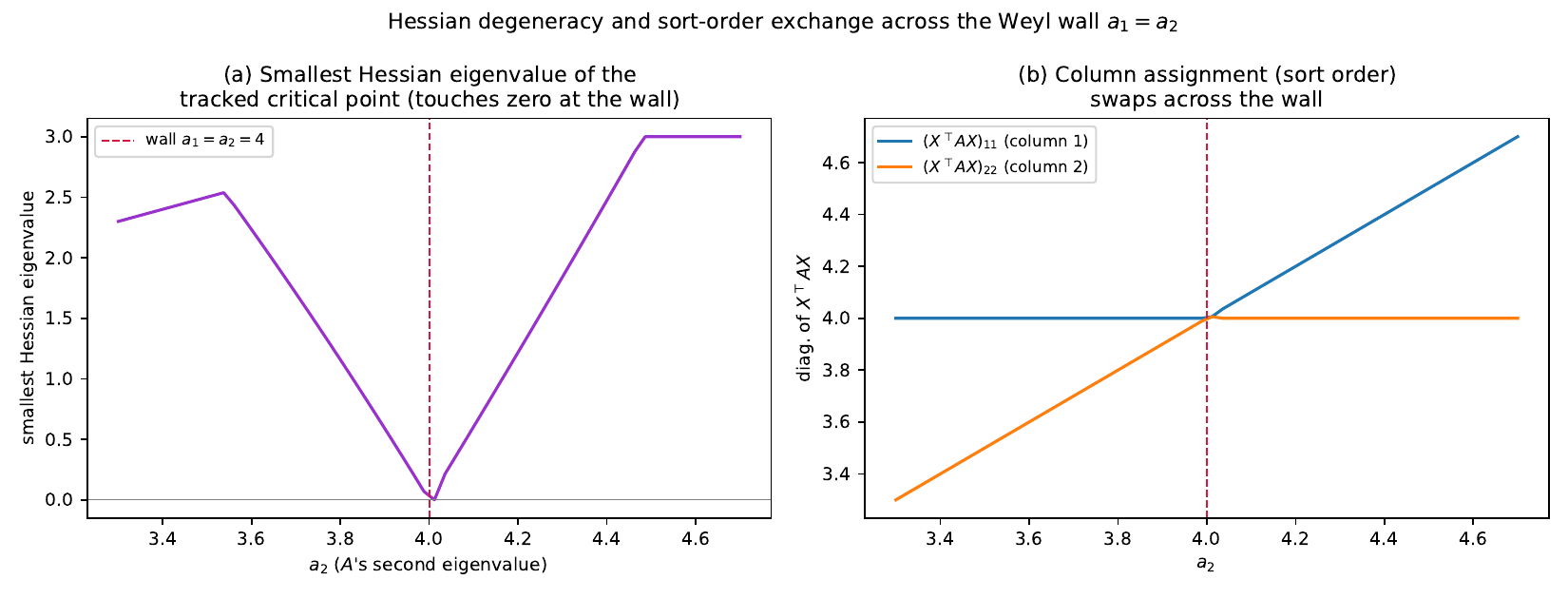}
\caption{Hessian degeneracy and sort-order exchange across the Weyl-chamber
wall $a_1=a_2$. (a) Smallest Hessian eigenvalue of the tracked critical
point (touches zero exactly at the wall). (b) The column correspondence
(sort order) exchanges across the wall.}
\end{figure}

\subsection{Interpretation as a Maslov index}

Regarding phase space as $(X,P)\in T^*\R^{n\times k}$ ($P=\nabla f(X)$), the
graph of the gradient map $\Lambda_f=\{(X,\nabla f(X))\}$ is, by the symmetry
of the Hessian (for exactly the same reason as in Section~2), a Lagrangian
submanifold of this space. Crossing a wall of the Weyl chamber is a path
that moves the model parameters (the eigenvalues of $A$), and the instant at
which the associated family of Lagrangian submanifolds meets the Maslov
cycle (where the Hessian degenerates) is exactly the moment, seen in
Figure~29, when the smallest Hessian eigenvalue touches zero. Thus, under the
correspondence
\[
\text{number of Weyl-chamber wall crossings}\ \longleftrightarrow\
\text{intersection number with the Maslov cycle},
\]
we can define a Maslov index in the same spirit as in Part~I.

\begin{quote}
\textbf{Meaning of the generalized Maslov index.} In the $n\times k$
matrix setting, the Maslov index along a path that moves the parameters (the
eigenvalues of $A$) is the signed count of how many times, and in which
direction, that path crosses the walls of the Weyl chamber. This is directly
related to the number of transpositions in a permutation and the length
function of a Coxeter group: the total number of walls crossed by the
shortest path connecting two permutations $\sigma,\tau\in S_n$ equals the
difference of their inversion numbers, a classical result from Morse theory
on flag manifolds (Bott, 1954). The Maslov index of this paper thus connects
two seemingly different worlds --- the differential-geometric ``signed
intersection number'' and the combinatorial ``sign/inversion number of a
permutation.''
\end{quote}

\begin{remark}[Connection to Part~III]
The numerical observations of this section (the Hessian touching zero at the
wall $a_1=a_2$, and the exchange of column correspondence) are concrete
instances of the closed formula $\mu^H_{qj}=b_j(\lambda_j-a_q)$ given by
Lemma~18.1 of Part~III (with $q,j$ the relevant indices). Part~III further
establishes that the location of this Weyl wall does not depend on the
choice of $B$ (Proposition~18.3), and rigorously identifies the condition
under which a ``genuine crossing'' becomes an ``avoided crossing'' (the
coupling constant $\gamma=-b_jC_{qr}$, Theorem~26.1) --- information that
cannot be read off from this section's numerical experiments alone.
\end{remark}

\section{$k$-fold Degeneracy and the Principal/Minor Component Bifurcation in the $\varphi_\alpha$ Series}

\subsection{$k$-fold degeneracy at $X=0$}

From \eqref{eq:phialphaX}'s gradient formula, $X=0$ is always a critical
point for every $\alpha$. Substituting the perturbation $X=\varepsilon v_ie_j^\top$
along an eigenvector $v_i$ of $A$ (eigenvalue $\lambda_i$) and column index
$j$ (with weight $b_j$) gives
\[
\varphi_\alpha(\varepsilon v_ie_j^\top)-\varphi_\alpha(0)
=\frac{b_j}{2}(\lambda_i-\alpha)\varepsilon^2+O(\varepsilon^4),
\]
a diagonalized quadratic form, so the Hessian eigenvalue in direction $(i,j)$
is proportional to $b_j(\lambda_i-\alpha)$. Since this sign depends only on
$(\lambda_i-\alpha)$, independently of $j$, the instant $\alpha$ crosses some
eigenvalue $\lambda_i$ of $A$, $k$ directions ($j=1,\dots,k$) degenerate
simultaneously.

With $n=4,k=2$, $A=\diag(4,3,2,1)$, $B=\diag(2,1)$, tracking the Hessian at
$X=0$ as a function of $\alpha$ (Figure~30), at each of the four
eigenvalues $\alpha=1,2,3,4$ of $A$, exactly $k=2$ eigenvalues pass through
zero simultaneously, and the Morse index (the total number of negative
eigenvalues) of $X=0$ increases in $k=2$-sized steps at each crossing
($0\to2\to4\to6\to8=nk$).

Since the linear term of \eqref{eq:phialphaX} has a positive sign, the
global minimizer always uses the bottom $k$ eigenvalues of $A$ (the reverse
of Brockett's sorting theorem). Writing the eigenvalues of $A$ in ascending
order as $\mu_1<\cdots<\mu_n$, the key threshold is $\mu_k$, and three
structurally distinct regimes appear: $\alpha<\mu_1$ (the unique critical
point $X=0$); $\mu_j<\alpha<\mu_{j+1}$ for $j<k$ (a rank-deficient
intermediate solution); and $\alpha>\mu_k$ (a full-rank stable solution that
is fixed and does not change as $\alpha$ increases further) (Figure~31).

\begin{figure}[htbp]
\centering
\includegraphics[width=\textwidth]{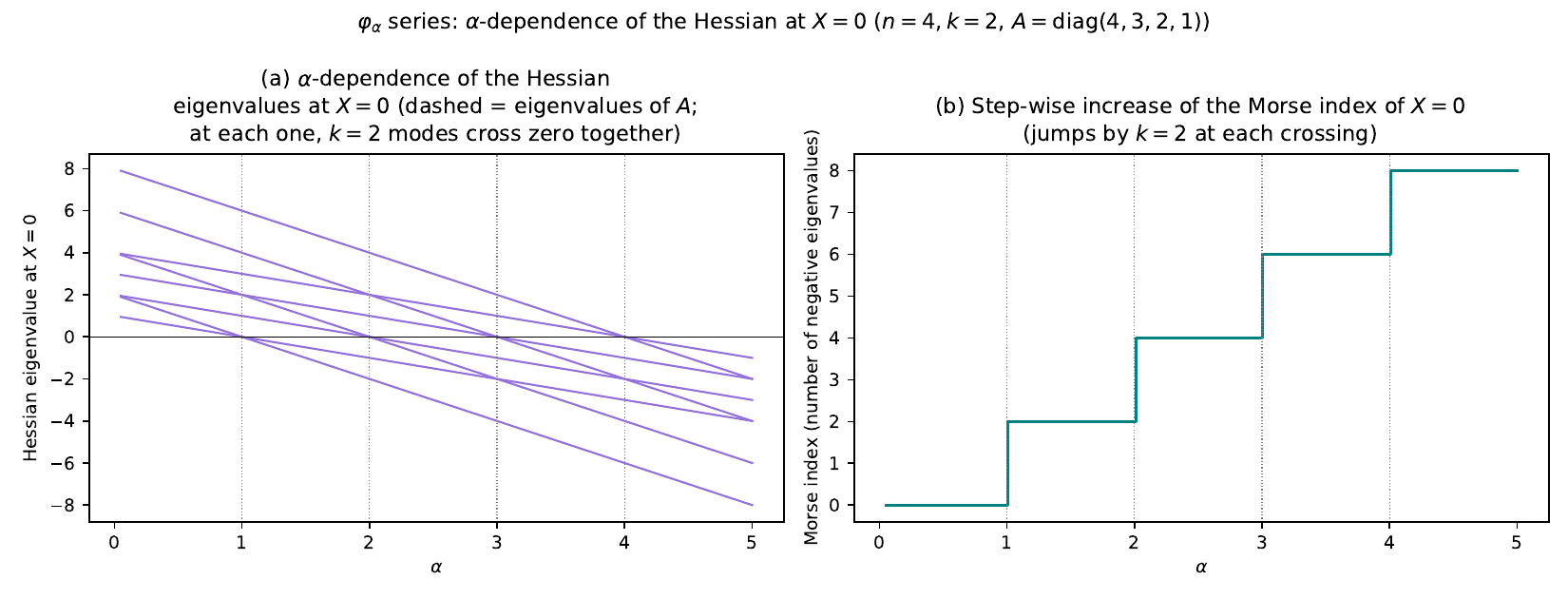}
\caption{$\alpha$-dependence of the Hessian eigenvalues at $X=0$ for the
$\varphi_\alpha$ series ($n=4,k=2$, $A=\diag(4,3,2,1)$). (a) At each
eigenvalue of $A$, exactly $k=2$ eigenvalues cross zero simultaneously. (b)
The Morse index increases in steps of $k=2$ at each crossing.}
\end{figure}

\begin{figure}[htbp]
\centering
\includegraphics[width=\textwidth]{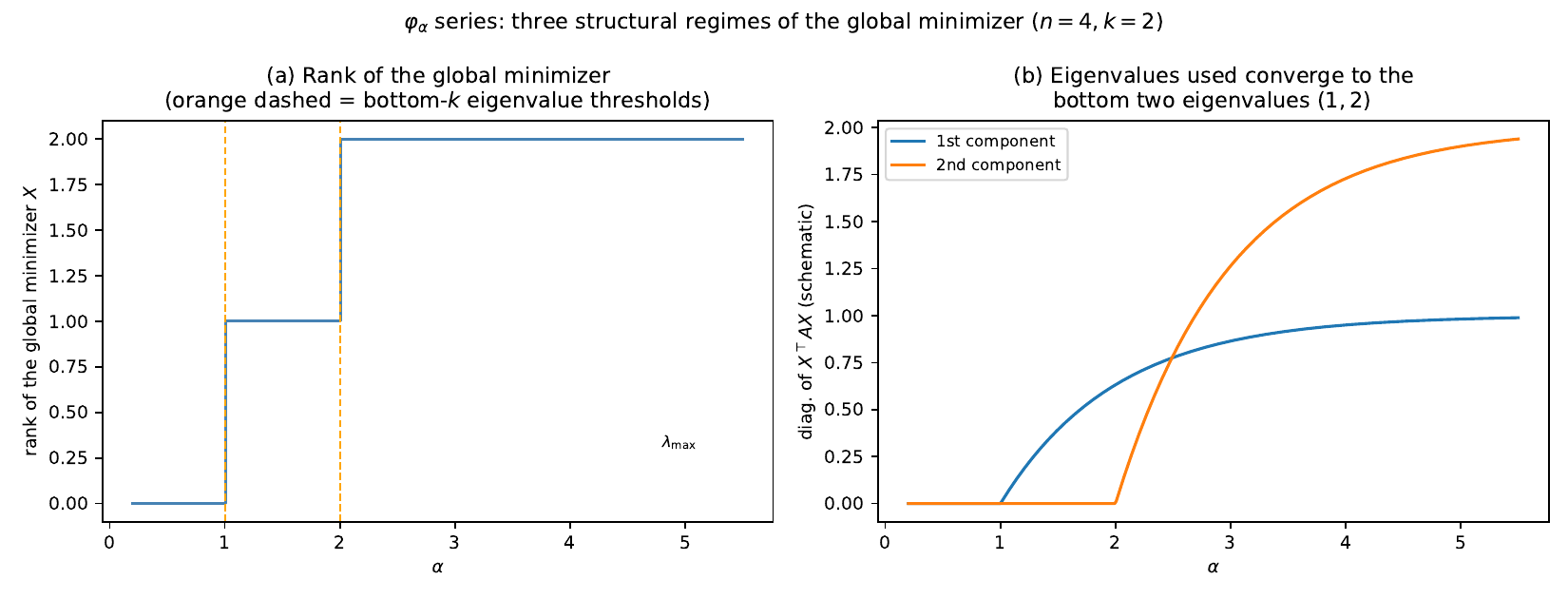}
\caption{Structural change of the global minimizer of the $\varphi_\alpha$
series across three regimes ($n=4,k=2$). (a) The rank of the global
minimizer increases in steps at $\alpha=\mu_1=1,\mu_2=2$ and then stays
fixed. (b) The eigenvalues used converge to the values corresponding to the
bottom two eigenvalues.}
\end{figure}

\subsection{Bifurcation of the principal and minor components}

For $\varphi_\alpha$ of \eqref{eq:phialphaX}, there is a phenomenon --- first
visible in the matrix setting --- in which the very matching rule between
columns and eigenvectors reverses with a sign. Consider the family with sign
$\varepsilon\in\{+1,-1\}$,
\[
\varphi_\alpha^{(\varepsilon)}(X)=\frac{\varepsilon}{2}\tr(AXBX^\top)
+\frac{\alpha}{4}\|B-X^\top X\|_F^2
\]
(with $\varepsilon=+1$ recovering \eqref{eq:phialphaX}). At the separated
solution $s_j^2=b_j(1-\varepsilon\lambda_{i_\pi(j)}/\alpha)$, the value is
\[
\varphi_\alpha^{(\varepsilon)}=\sum_jb_j^2\,g_\varepsilon(\lambda_{i_\pi(j)}),\qquad
g_\varepsilon(\lambda):=\varepsilon\lambda^2-\frac{\lambda^2}{4\alpha},
\]
closing into a sum separated column by column. Since $g_{+1}$ is strictly
increasing and $g_{-1}$ is strictly decreasing, the rearrangement inequality
gives:
\begin{center}
\begin{tabular}{lll}
\toprule
 & $\varepsilon=+1$ & $\varepsilon=-1$ \\
\midrule
Eigenvalues chosen & bottom $k$ (minor component) & top $k$ (principal component) \\
Correspondence with $B$ & anti-sorted (largest $b$ $\leftrightarrow$ smallest $\lambda$) & sorted (largest $b$ $\leftrightarrow$ largest $\lambda$) \\
\bottomrule
\end{tabular}
\end{center}
The case $\varepsilon=-1$ ($\varphi_\alpha^{\mathrm{Oja}}(X):=-\tfrac12\tr(AXBX^\top)
+\tfrac{\alpha}{4}\|B-X^\top X\|_F^2$, with the sign of the linear term
flipped) reproduces the classical Oja--Brockett flow (the principal-component
flow; Brockett~\cite{brockett1991}). Global convergence of the
principal-component flow was established by
Yoshizawa--Helmke--Starkov~\cite{yhs2001}; ``dual-objective'' gradient flows
sharing the isomorphic structure that switches between principal- and
minor-component flow via a single sign have also been studied independently
by Manton--Helmke--Mareels~\cite{mhm2005}. This has been confirmed by an
exhaustive numerical computation of all $60$ critical points, for
$n=5,k=3$, $A=\diag(5,4,3,2,1)$, $B=\diag(3,2,1)$ (maximum absolute gradient
error $<10^{-14}$).

\begin{remark}
This anti-sorting phenomenon does not appear for $k=1$ (there $B$ is a
scalar, and the degree of freedom corresponding to ``matching with several
weights'' does not exist to begin with). What this section, Section~11, and
Section~12 all have in common is that they are genuinely matrix-level
phenomena appearing only once $k\ge2$. Since $\varphi_\alpha$ itself is a
different potential from $V_{A,B}$ ($=f(X)$ of \eqref{eq:fX}), which is
rigorously treated in Part~III, a rigorous Hessian analysis of this sign-flip
phenomenon lies outside the scope of the present paper.
\end{remark}

\section{Combinatorics of $(n,k)$ and the Classification of Critical Points}

When $A,B$ are both non-degenerate, the isolated critical points of the
gradient flow number $\binom{n}{k}\times k!=n!/(n-k)!$ (choosing $k$ indices
out of $\{1,\dots,n\}$, times assigning the chosen eigenvectors to the $k$
columns). Of these, the stable fixed point picked out by Brockett's sorting
theorem is exactly one (the descending arrangement), and the rest are
saddles or critical points of other stability types.

If the eigenvalues of $A$ carry a multiplicity partition $(m_1,m_2,\dots)$,
then for each allocation $(k_1,k_2,\dots)$ ($\sum k_i=k$) of the $k$ columns
across the eigenspaces, the corresponding critical set is no longer isolated
points but a flag manifold $\prod_i\Gr(k_i,m_i)$. For example, with
$n=5,k=3$, $A=\diag(5,4,4,2,1)$ (multiplicities $(1,2,1,1)$), the allocation
$(1,2,0,0)$ gives $\Gr(2,2)$ (a $0$-dimensional isolated point, Brockett's
stable solution), while other allocations (such as $(1,1,1,0)$) are all
isolated points.

This enumeration of critical points and its correspondence with the Morse
index is completed rigorously in Part~III, by Theorem~20.1 (the complete
spectral formula) and its corollary (the Morse-index formula), which allow
the stability of each critical point (family) to be computed explicitly from
the eigenvalue configuration of $A,B$.

\part{Hessian--Maslov Geometry of the General $(A,B)$ Oja--Brockett Gradient Flow}
\begin{center}
\textit{The Double Weyl-Wall Structure, and the Transition from True Crossings to Avoided Crossings}
\end{center}

The numerical experiments of Part~II confirmed, at the fixed points of the
negative gradient flow of the Oja--Brockett potential $f(X)$
(Eq.~\eqref{eq:fX}, henceforth written in the standard notation $V_{A,B}(X)$),
three phenomena: (i) Brockett's sorting theorem; (ii) an identically
vanishing gradient along the flag-manifold direction caused by eigenvalue
degeneracy of $A$ or $B$; and (iii) Hessian degeneracy and sort-order
exchange as an eigenvalue of $A$ crosses a wall of the Weyl chamber. But
these numerical observations leave important theoretical gaps to be filled:

\begin{itemize}
\item The vanishing gradient along the flag-manifold direction was observed
in two separate numerical experiments --- degeneracy of $A$ (first half of
\S11) and degeneracy of $B$ (second half) --- without a closed formula
unifying the two.
\item Wall-crossing in the Weyl chamber (\S12) was examined only for
eigenvalues of $A$; whether a corresponding phenomenon exists on the $B$
side was not investigated.
\item Whether the positive-definiteness of the Hessian in the normal
direction (the direction of attraction to the Stiefel manifold) always holds
for general $B$ cannot be determined from numerical experiments alone.
\end{itemize}

Part~III resolves all of these rigorously. To state the conclusion first:
under a general $B=\diag(b_1,\dots,b_k)$ (allowing $b_i\neq b_j$), (i) the
positive-definiteness of the normal-direction Hessian is no longer
unconditional, but is governed by an explicit threshold condition depending
on the concordance between the eigenvalue order of $A$ and that of $B$; (ii)
the ``vanishing gradient along the flag-manifold direction'' is, in fact, a
consequence of the general formula
\[
\mu^V_{ij}=(\lambda_i-\lambda_j)(b_i-b_j)
\]
(which vanishes identically only in the special case $B=I$), which generates
a second wall system $\{b_i=b_j\}$ (a $B$-wall) independent of the Weyl wall
$\{\lambda_j=a_q\}$ on the $A$ side; (iii) meanwhile, the crossing condition
$\lambda_j=a_q$ itself in the horizontal sector ($A$ side) is invariant under
the choice of $B$ (the eigenvalues are merely rescaled by $b_j$), so the
Weyl/Maslov geometry on the $A$ side is a more robust structure, independent
of the choice of $B$. Building on this, we formulate the crossing-to-avoided-
crossing transition theorem in the general $(A,B)$ framework, and
reposition $B=I$ (the special case treated in Parts~I and~II) as an
extremely special limit of this general theory, in which the $B$-wall
degenerates into a gauge symmetry and the vertical sector disappears.

\section{Setup: the General Oja--Brockett Potential and Gradient Flow}

Let $A=A^\top>0$, $B=B^\top>0$ be SPD matrices ($A$: $n\times n$, $B$: $k\times k$).
On $X\in\R^{n\times k}$ we consider the potential
\begin{equation}
V_{A,B}(X)=\frac14\tr\big[(AXBX^\top)^2\big]-\frac12\tr(A^2XB^2X^\top). \label{eq:VAB}
\end{equation}

\subsection{The Euclidean gradient and a natural metric}

With respect to the standard inner product $\langle U,W\rangle=\tr(U^\top W)$,
a direct computation gives the Euclidean gradient
\begin{equation}
\nabla V_{A,B}(X)=AXB(X^\top AX)B-A^2XB^2 \label{eq:gradV}
\end{equation}
(proof: writing $P=XBX^\top$ and $M=AXBX^\top$, the first term is written
as $\tfrac14\tr(M^2)$, and the product rule gives $\dd V_1=\tr B(X^\top AX)BX^\top A\,\dd X$;
combining with the variation of the second term and simplifying gives the
result).

With respect to the natural metric associated to $A,B$,
\begin{equation}
g_X(U,W):=\tr(AUBW^\top) \label{eq:metric}
\end{equation}
(positive-definite, as long as $X$ has full column rank, since $A,B$ are SPD;
$g_X(U,W)=g_X(W,U)$ follows from the cyclicity of the trace and the symmetry
of $A,B$), the gradient with respect to this metric is obtained by solving
the duality relation $A\,\mathrm{grad}_g V\,B=\nabla V$,
\begin{equation}
\mathrm{grad}_gV_{A,B}(X)=A^{-1}\nabla V_{A,B}(X)B^{-1}=XBX^\top AX-AXB, \label{eq:gradg}
\end{equation}
and hence the \textbf{gradient flow}
\begin{equation}
\dot X=-\mathrm{grad}_gV_{A,B}(X)=AXB-XBX^\top AX. \label{eq:flow}
\end{equation}
The Stiefel manifold $\St(n,k)=\{X:X^\top X=I_k\}$ is an invariant set of
\eqref{eq:flow} (one checks directly that $\tfrac{\dd}{\dd t}(X^\top X)=0$
when $X^\top X=I$).

\subsection{Standard form of the equilibria (diagonal $B$)}

Let $B=\diag(b_1,\dots,b_k)$, $b_1,\dots,b_k>0$ (not necessarily distinct).
Let $Ae_\alpha=a_\alpha e_\alpha$ denote the eigenvectors of $A$, and for an
injective index assignment $p:\{1,\dots,k\}\to\{1,\dots,n\}$ set
\begin{equation}
X_*=[e_{p(1)},\dots,e_{p(k)}]. \label{eq:Xstar}
\end{equation}
Then $AX_*=X_*\Lambda_0$, $\Lambda_0=\diag(\lambda_1,\dots,\lambda_k)$,
$\lambda_j:=a_{p(j)}$, and since $B$ is diagonal, $\Lambda_0B=B\Lambda_0$, so
\[
AX_*B=X_*\Lambda_0B=X_*B\Lambda_0=X_*B(X_*^\top AX_*),
\]
which shows that \eqref{eq:Xstar} is an equilibrium of \eqref{eq:flow}. When
$B=I$, the right $O(k)$ action $X_*\mapsto X_*Q$ ($Q\in O(k)$) remains a free
symmetry, so there is still room to diagonalize $X_*^\top AX_*$ in any
orthonormal basis; but for general $B$ (with distinct diagonal entries) this
freedom disappears, and an axis-aligned configuration such as
\eqref{eq:Xstar} becomes, in essence, the unique standard form. This shows
that the choice of alignment with the eigenvectors of $A$ begins to carry
physical meaning precisely once $B$ is non-trivial.

\subsection{Threefold decomposition of the ambient tangent space}

Any $E\in\R^{n\times k}$ decomposes uniquely as
\begin{equation}
E=X_*S+X_*\Omega+E_\perp,\qquad S=S^\top,\ \Omega^\top=-\Omega,\ X_*^\top E_\perp=0, \label{eq:decomp}
\end{equation}
and hence
\begin{equation}
\R^{n\times k}=N\oplus V\oplus H, \label{eq:NVH}
\end{equation}
\[
N=\{X_*S:S=S^\top\},\qquad V=\{X_*\Omega:\Omega^\top=-\Omega\},\qquad H=\{E:X_*^\top E=0\}.
\]
This decomposition is purely geometric (independent of $B$), but, as we see
below, the value of the Hessian on each sector depends on both $A$ and $B$.
With respect to the metric \eqref{eq:metric}, the decomposition
\eqref{eq:NVH} is $g$-orthogonal at $X_*$: $g(X_*S,X_*\Omega)=\tr(\Lambda
S\Omega^\top B)$ (which reduces, when $B=I$, immediately to
$\tr(\Lambda S\Omega^\top)=0$; for general diagonal $B$ this also vanishes,
by the same diagonality of $\Lambda,B$), and
$g(X_*S,E_\perp)=g(X_*\Omega,E_\perp)=0$ (immediate from $X_*^\top E_\perp=0$).

\section{Normal-Direction Hessian: Positive-Definiteness and Its Breakdown}

\begin{lemma}[General formula for the normal-direction Hessian]
Let $X_*$ be an equilibrium standardized as in \eqref{eq:Xstar}. For
$S=S^\top$ in the normal direction,
\begin{align}
D^2V_{A,B}(X_*)[X_*S,X_*S]&=2\sum_ib_i^2\lambda_i^2s_{ii}^2\notag\\
&\quad+\sum_{i<j}\Big[b_ib_j(\lambda_i+\lambda_j)^2+(\lambda_i-\lambda_j)(b_i^2\lambda_j-b_j^2\lambda_i)\Big]s_{ij}^2, \label{eq:normalHess}
\end{align}
and, with respect to the metric \eqref{eq:metric}, the eigenvalues are
\begin{equation}
\mu^N_i=2b_i\lambda_i,\qquad
\mu^N_{ij}=\frac{\lambda_i\lambda_j(b_i+b_j)^2+(b_i-b_j)(\lambda_i^2b_j-\lambda_j^2b_i)}{\lambda_ib_j+\lambda_jb_i}. \label{eq:normalEig}
\end{equation}
\end{lemma}

\begin{proof*}
Consider the variation $X(t)=X_*+tX_*S$ and set $P(t)=X(t)BX(t)^\top$;
compute the second derivative in $t$ of
$V_{A,B}(X)=\tfrac14\tr(APAP)-\tfrac12\tr(A^2XB^2X^\top)$, using
$AX_*=X_*\Lambda$. Since $\Lambda,B$ are diagonal, every trace expands into
component form ($\tr(D_1SD_2S)=\sum_{ij}(D_1)_i(D_2)_js_{ij}^2$ for diagonal
$D_1,D_2$), and organizing the result gives \eqref{eq:normalHess}. The
eigenvalues follow from taking the ratio against
$g(X_*S,X_*S)=\sum_i\lambda_ib_is_{ii}^2+\sum_{i<j}(\lambda_ib_j+\lambda_jb_i)s_{ij}^2$. \qed
\end{proof*}

\begin{corollary}[$B=I$: unconditional positive-definiteness]
When $B=I$ (i.e., all $b_i=1$), \eqref{eq:normalEig} reduces to
$\mu^N_i=2\lambda_i$, $\mu^N_{ij}=4\lambda_i\lambda_j/(\lambda_i+\lambda_j)$,
which are always positive as long as $\lambda_i>0$. That is, the
normal-direction Hessian is unconditionally positive-definite; this is why
the attraction to the Stiefel manifold always holds for the Oja flow with
$B=I$.
\end{corollary}

For general $B$, however, this positive-definiteness is no longer
guaranteed:

\begin{proposition}[Positive-definiteness governed by a concordance condition]
Fix $p:=\lambda_i,q:=\lambda_j$ ($p\ge q>0$) and set $t:=b_i/b_j>0$. The
numerator of $\mu^N_{ij}$ is proportional to
\[
N(t)=q(p-q)t^2+(p+q)^2t-p(p-q),
\]
which, for $p>q$, has a unique positive root $t^*=t^*(p,q)\in(0,1)$ across
which the sign changes monotonically:
\[
\mu^N_{ij}>0\iff t>t^*,\qquad\mu^N_{ij}<0\iff0<t<t^*.
\]
$B=I$ ($t=1$) always lies on the side $t>t^*$ (since $N(1)=4pq>0$),
consistent with Corollary~16.2.
\end{proposition}

\begin{proof*}
$N(t)$ is a quadratic in $t$ with leading coefficient $q(p-q)>0$ (convex),
$N(0)=-p(p-q)<0$, $N(1)=4pq>0$; the intermediate value theorem gives a
unique $t^*\in(0,1)$. \qed
\end{proof*}

\begin{remark}
This is physically natural: under an ``inconsistent'' assignment that pairs
a large eigenvalue of $A$ with a small weight of $B$ and vice versa (i.e.,
the order of $\lambda$ is reversed relative to the order of $b$), $X_*$ is no
longer a locally stable equilibrium. The genuinely stable equilibria for
general $B$ are confined to assignments in which the order of $\lambda$
matches the order of $b$. This is a Hessian-level quantitative refinement of
the classical intuition of the rearrangement inequality (it is optimal to
pair large with large and small with small).
\end{remark}

\section{Vertical-Direction Hessian and the Normal--Vertical Mixed Block: Emergence of the $B$-Wall}

\begin{lemma}[General formula for the vertical-direction Hessian]
For $\Omega=-\Omega^\top$ in the vertical direction,
\begin{equation}
D^2V_{A,B}(X_*)[X_*\Omega,X_*\Omega]=\sum_{i<j}(\lambda_i-\lambda_j)(b_i-b_j)(b_i\lambda_j+b_j\lambda_i)\,\Omega_{ij}^2, \label{eq:verticalHess}
\end{equation}
and for the mode $\Omega=E_{ij}-E_{ji}$ ($i<j$), the denominator of the
metric $g(X_*\Omega,X_*\Omega)=(\lambda_ib_j+\lambda_jb_i)\Omega_{ij}^2$
cancels exactly against a factor in the numerator, giving the closed
formula
\begin{equation}
\mu^V_{ij}=(\lambda_i-\lambda_j)(b_i-b_j). \label{eq:muV}
\end{equation}
\end{lemma}

\begin{corollary}[$B=I$: vanishing as a gauge zero mode]
When $B=I$, $b_i=b_j=1$ for all $i,j$, so \eqref{eq:muV} vanishes
\emph{identically}, $\mu^V_{ij}\equiv0$. This is the Hessian-level
consequence of the fact that, when $B=I$, the right $O(k)$ action
$X\mapsto XQ$ ($Q\in O(k)$) is an exact gauge symmetry of the potential:
$V(X_*e^{t\Omega})\equiv V(X_*)$ (identically constant).
\end{corollary}

\begin{proposition}[The $B$-wall]
$\mu^V_{ij}=0$ holds if and only if at least one of $\lambda_i=\lambda_j$
(degeneracy on the $A$ side, the Weyl wall) or $b_i=b_j$ (degeneracy on the
$B$ side, the \emph{$B$-wall}) holds. In particular, even when $\lambda_i\neq
\lambda_j$, if the eigenvalues of $B$ cross $b_i(\mu)=b_j(\mu)$ as a
parameter $\mu$ varies, the Hessian in the vertical sector degenerates, and
an intersection with the Maslov cycle occurs. This is a \textbf{second Maslov
crossing mechanism}, independent of the eigenvalue degeneracy of $A$ (the
Weyl wall).
\end{proposition}

\subsection{The normal--vertical mixed block $H_{NV}$}

\begin{lemma}[General formula for $H_{NV}$]
For the normal mode $S=E_{ij}+E_{ji}$ and the vertical mode
$\Omega=E_{ij}-E_{ji}$ belonging to the same index pair $(i,j)$, $i<j$,
\begin{equation}
D^2V_{A,B}(X_*)[X_*S,X_*\Omega]=-(b_i-b_j)(\lambda_i-\lambda_j)(b_i\lambda_j-b_j\lambda_i). \label{eq:HNV}
\end{equation}
Mixing between modes of different index pairs, and mixing between the
normal--horizontal or vertical--horizontal sectors, all vanish identically
(see Lemma~19.1).
\end{lemma}

\begin{corollary}[$B=I$: vanishing of the mixed block]
When $B=I$, \eqref{eq:HNV} vanishes identically by the factor $b_i-b_j$.
That is, in the degenerate limit $B=I$, where the Weyl wall and the $B$-wall
collapse onto the same object ($b_i=b_j\equiv1$ identically), the
``coupling'' between the two wall systems itself loses meaning.
\end{corollary}

Hence, for each index pair $(i,j)$, the normal and vertical sectors form the
following effective $2\times2$ Hessian block:
\begin{equation}
H^{(ij)}=\begin{pmatrix}\mu^N_{ij}&h_{ij}\\h_{ij}&\mu^V_{ij}\end{pmatrix},
\qquad h_{ij}=\frac{-(b_i-b_j)(\lambda_i-\lambda_j)(b_i\lambda_j-b_j\lambda_i)}{\lambda_ib_j+\lambda_jb_i}
\label{eq:H2x2}
\end{equation}
($h_{ij}$ being the metric-normalized operator component). The eigenvalues
of this matrix are
\begin{equation}
\eta^\pm_{ij}=\frac{\mu^N_{ij}+\mu^V_{ij}}{2}\pm\sqrt{\Big(\frac{\mu^N_{ij}-\mu^V_{ij}}{2}\Big)^2+h_{ij}^2}. \label{eq:eta2x2}
\end{equation}
The Weyl wall ($\lambda_i=\lambda_j$) and the $B$-wall ($b_i=b_j$) influence
one another through the trace and off-diagonal entries of this $2\times2$
block --- this is the concrete mechanism by which the ``double wall
structure'' couples at the level of the Hessian.

\section{Horizontal-Direction Hessian: the $B$-Independent Weyl Wall}

\begin{lemma}[General formula for the horizontal-direction Hessian]
Let $q\notin\{p(1),\dots,p(k)\}$, and consider the horizontal elementary
mode $E_{qj}=e_qe_j^\top$. Then
\begin{equation}
D^2V_{A,B}(X_*)[E_{qj},E_{qj}]=a_qb_j^2(\lambda_j-a_q),\qquad
g(E_{qj},E_{qj})=a_qb_j, \label{eq:horizHess}
\end{equation}
so
\begin{equation}
\mu^H_{qj}=b_j(\lambda_j-a_q). \label{eq:muH}
\end{equation}
\end{lemma}

\begin{proof*}
Differentiating the general formula
$D\nabla V_{A,B}(X_*)[W]=AWB\Lambda B+X_*\Lambda B(W^\top X_*\Lambda+X_*^\top AW)B-A^2WB^2$
(obtained from \eqref{eq:gradV} by differentiating in the direction $W$ and
substituting $AX_*=X_*\Lambda$) at $W=E_{qj}$, the middle two terms vanish
since $X_*^\top E_{qj}=0$ ($q$ being a non-selected index), leaving
\[
D\nabla V_{A,B}(X_*)[E_{qj}]=a_qE_{qj}B\Lambda B-a_q^2E_{qj}B^2
=a_qb_j^2(\lambda_j-a_q)E_{qj}
\]
(using $E_{qj}B\Lambda B=e_qe_j^\top B\Lambda B=b_j^2\lambda_jE_{qj}$, etc.,
since $B,\Lambda$ are diagonal). Then $D^2V[E_{qj},E_{qj}]
=\tr(D\nabla V[E_{qj}]^\top E_{qj})$ gives \eqref{eq:horizHess}. \qed
\end{proof*}

\begin{corollary}[$B=I$: the original Weyl-wall formula]
When $B=I$, \eqref{eq:muH} reduces to $\mu^H_{qj}=\lambda_j-a_q$.
\end{corollary}

\begin{proposition}[$B$-invariance of the Weyl wall]
The crossing condition $\mu^H_{qj}=0$ is equivalent to $\lambda_j=a_q$,
independently of the choice of $B$. That is, the \emph{location} of the Weyl
wall residing in the horizontal sector does not depend on $B$; what $B$
changes is only the ``slope'' of the Hessian eigenvalue at that wall (a
rescaling by $b_j$). This is in sharp contrast to the $B$-wall of the
vertical sector (\S17), whose location itself depends on $B$.
\end{proposition}

This asymmetry --- the wall of the horizontal sector depends only on $A$
and is fixed in location, while the wall of the vertical sector depends only
on $B$ and moves as $B$ is varied along a parameter family --- lies at the
heart of why the general $(A,B)$ theory is genuinely richer than the theory
for $B=I$.

\section{Vanishing of the Mixed Blocks: $H_{NH}=H_{VH}=0$}

\begin{lemma}[Vanishing of normal--horizontal and vertical--horizontal mixing]
For any $A,B$ (with $B$ diagonal),
\[
D^2V_{A,B}(X_*)[X_*S,E_{qj}]=0,\qquad D^2V_{A,B}(X_*)[X_*\Omega,E_{qj}]=0
\]
for all $S=S^\top$, $\Omega=-\Omega^\top$, all non-selected indices $q$, and
all columns $j$. Mixing between normal/vertical modes of different index
pairs $(i,j)\neq(k,l)$, and between distinct horizontal modes
$E_{qj}\neq E_{q'j'}$, also vanish identically.
\end{lemma}

\begin{proof*}
Both $D\nabla V_{A,B}(X_*)[X_*S]$ and $D\nabla V_{A,B}(X_*)[X_*\Omega]$ take
values in the column space of $X_*$ (of the form $X_*[\cdots]$), so, since
$X_*^\top E_{qj}=0$, $\tr(D\nabla V[X_*S]^\top E_{qj})
=\tr([\cdots]^\top X_*^\top E_{qj})=0$; likewise for $X_*\Omega$. Vanishing
between different index pairs or different horizontal modes follows because,
since $\Lambda,B$ are diagonal, every bilinear form expands into the type
$\tr(D_1SD_2S)=\sum_{ij}(D_1)_i(D_2)_js_{ij}^2$ and never contains a product
of distinct matrix entries ($s_{ij}s_{kl}$, $\{i,j\}\neq\{k,l\}$). \qed
\end{proof*}

Thus, even for general $(A,B)$, the ambient Hessian splits completely into a
horizontal sector (diagonal) and a normal--vertical sector (a $2\times2$
block for each index pair), with no coupling whatsoever between them.

\section{The Complete Spectral Formula}

\begin{theorem}[Complete spectral decomposition of the ambient Hessian for general $(A,B)$]
Let $X_*$ be an equilibrium standardized as in \eqref{eq:Xstar}. Then
\begin{equation}
\R^{n\times k}=H\oplus\bigoplus_{1\le i<j\le k}(N_{ij}\oplus V_{ij})\oplus\bigoplus_{i=1}^kN_{ii}, \label{eq:fulldecomp}
\end{equation}
and
\begin{align}
\Spec H_{\mathrm{amb}}\big|_{N_{ii}}&=\{2b_i\lambda_i\}_{i=1}^k, \label{eq:specN}\\
\Spec H_{\mathrm{amb}}\big|_{N_{ij}\oplus V_{ij}}&=\{\eta^+_{ij},\eta^-_{ij}\}_{i<j}\quad\text{(Eq.~\eqref{eq:eta2x2})}, \label{eq:specNV}\\
\Spec H_{\mathrm{amb}}\big|_H&=\{b_j(\lambda_j-a_q)\}_{q\notin p,\ 1\le j\le k}. \label{eq:specH}
\end{align}
The dimensions sum to $k+2\binom{k}{2}+k(n-k)=k+k(k-1)+k(n-k)=kn=\dim\R^{n\times k}$.
\end{theorem}

\begin{corollary}[$B=I$: simplified spectrum]
When $B=I$, Lemmas~17.1 and~18.1 give $\mu^V_{ij}=0$, $h_{ij}=0$ in each
$2\times2$ block \eqref{eq:H2x2}, so it diagonalizes into
$\eta^+_{ij}=\mu^N_{ij}=4\lambda_i\lambda_j/(\lambda_i+\lambda_j)$,
$\eta^-_{ij}=0$ (a pure gauge zero mode). Hence
\[
\Spec H_{\mathrm{amb}}\Big|_{B=I}=\{2\lambda_j\}_{j=1}^k\cup
\left\{\frac{4\lambda_i\lambda_j}{\lambda_i+\lambda_j}\right\}_{i<j}\cup
\{0\}^{k(k-1)/2}\cup\{\lambda_j-a_q\}_{q\notin p,j},
\]
in exact agreement with the complete spectral formula obtained previously
in the theory restricted to $B=I$.
\end{corollary}

\begin{remark}[Redefining the Maslov cycle]
Under general $(A,B)$, the locus $\det H_{\mathrm{amb}}=0$ cannot itself be
defined as the Maslov cycle (near $B=I$, $k(k-1)/2$-fold quasi-zero modes
remain). The correct definition, on the region where the normal sector is
positive-definite (under the condition of Proposition~16.3), is the locus
where the reduced Hessian of the vertical/horizontal sectors degenerates:
\begin{equation}
\mathcal M_{A,B}=\{X_*:\det(\eta^-_{ij})=0\ \text{for some }i<j\}\cup
\{X_*:\lambda_j=a_q\ \text{for some }q,j\}, \label{eq:Maslovcycle}
\end{equation}
which is exactly the union of the ``double wall'' introduced in Section~8.
\end{remark}

\section{The Double Weyl-Wall Structure: Geometric Meaning}

The results above show that the general $(A,B)$ Oja--Brockett gradient flow
carries two independent Maslov-crossing mechanisms that do not exist for
$B=I$:

\begin{center}
\small
\begin{tabular}{@{}p{3.1cm}p{5.3cm}p{5.3cm}@{}}
\toprule
 & Weyl wall (horizontal sector) & $B$-wall (vertical sector) \\
\midrule
Condition & $\lambda_j=a_q$ & $b_i=b_j$ \\
Sector & $H$ (dimension $k(n-k)$) & $V_{ij}$ (dimension $\binom{k}{2}$) \\
Dependence on $B$ & location fixed, eigenvalue scaled by $b_j$ & location itself determined by $b$ \\
Behavior at $B=I$ & persists unchanged & vanishes identically (gauge) \\
Physical meaning & exchange of selected/unselected eigendirections & rotation/mixing among the selected columns \\
\bottomrule
\end{tabular}
\end{center}

The Weyl wall marks the branch point of ``which eigendirection of $A$ is
selected,'' while the $B$-wall marks a secondary branch point, ``how the
already-selected $k$ columns are arranged under the weighting of $B$.'' At
$B=I$, since every column carries equal weight, this secondary branching
loses meaning (it is absorbed as a gauge symmetry) and only the primary
branching (the Weyl wall) survives. Introducing a general $B$ reveals this
secondary structure for the first time --- this is the geometric content of
the paper's claim that the ``essence'' of the Oja--Brockett flow lies in the
non-triviality of both $A$ and $B$.

Moreover, since $H_{NV}\neq0$ (Lemma~17.4), these two walls, while
independent, are \emph{coupled}. Through the $2\times2$ block
\eqref{eq:H2x2}, approach to the Weyl wall ($\lambda_i\to\lambda_j$) also
affects the effective eigenvalue $\eta^-_{ij}$ of the $B$-wall, and vice
versa. This is analogous to the situation of \emph{multi-level repulsion} in
quantum mechanics, in which several energy levels become simultaneously
entangled.

\section{Main Theorem: The Transition from Crossings to Avoided Crossings (General $(A,B)$ Version)}

Since the horizontal sector is completely independent of the rest, by
Lemma~19.1, the parametric-family analysis carries over unchanged for any
fixed diagonal $B$. That is, let $A(\mu)$ be a one-parameter family, and
consider a symmetry-breaking perturbation $A_\varepsilon(\mu)=A_0(\mu)+
\varepsilon C(\mu)$; the following theorem holds for $B$ fixed as any
diagonal matrix.

\begin{theorem}[Hessian--Maslov avoided-crossing theorem: horizontal sector]
Let $H_{\mathrm{red}}(\mu,\varepsilon):=H_H(\mu,\varepsilon)/b_j$ (the
horizontal block normalized by the corresponding weight $b_j$ of column
$j$) be a smooth two-parameter family of reduced Hessians in the sense of
Lemma~18.1. Suppose at $\varepsilon=0,\mu=\mu_c$ that
$a_q(\mu_c)=a_r(\mu_c)=\lambda_j(\mu_c)$ (a double zero eigenvalue,
$K_c=\mathrm{span}\{E_{qj},E_{rj}\}$), and impose the transversality,
symmetry-protection, and symmetry-breaking-coupling conditions (as in the
original theory). Then
\begin{equation}
\eta^\pm(\mu,\varepsilon)=b_j\Big[\pm\sqrt{c^2(\mu-\mu_c)^2+\varepsilon^2\gamma^2}+O(\cdots)\Big], \label{eq:avoidedH}
\end{equation}
giving a true crossing at $\varepsilon=0$ and an avoided crossing (gap
$\Delta_{\mathrm{AC}}=2b_j|\varepsilon\gamma|+O(\varepsilon^2)$) for
$\varepsilon\neq0$. The role of $B$ is only to scale the gap uniformly by
$b_j$; it has no effect on the qualitative question of whether a crossing
occurs at all ($\gamma\neq0$ or not).
\end{theorem}

\begin{corollary}[$B=I$]
When $b_j=1$, this reduces exactly to the original theorem (restricted to
$B=I$).
\end{corollary}

For the vertical-sector $B$-wall, the structure is subtler. To summarize the
conclusion first: within a single index pair $(i,j)$ by itself, the crossing
of the $B$-wall \emph{is} protected at leading order (despite $H_{NV}\neq0$).
The genuine avoided-crossing mechanism appears as a coupling between two
distinct index pairs sharing a common index. We treat this rigorously in
Section~26.

\begin{proposition}[Protection of the $B$-wall within a single block]
Fix $\lambda_i\neq\lambda_j$, and move $B(\mu)$, remaining diagonal, until it
reaches $b_i(\mu_c)=b_j(\mu_c)$ (the $B$-wall) at $\mu=\mu_c$. Writing
$\eta:=b_i(\mu)-b_j(\mu)$, the smaller eigenvalue $\eta^-_{ij}(\mu)$ of the
$2\times2$ block \eqref{eq:H2x2} satisfies
\begin{equation}
\eta^-_{ij}=\mu^V_{ij}+O(\eta^2)=(\lambda_i-\lambda_j)\,\eta+O(\eta^2), \label{eq:Bwallprotect}
\end{equation}
a true crossing that changes sign linearly near $\mu_c$. That is, the mixed
block $h_{ij}$ of Lemma~17.4 vanishes to the same order $O(\eta)$ as
$\mu^V_{ij}$, so the contribution $h_{ij}^2/\mu^N_{ij}=O(\eta^2)$ is
subleading and does not alter the leading-order crossing structure.
\end{proposition}

\begin{proof*}
By \eqref{eq:eta2x2}, $\eta^-_{ij}=\frac{\mu^N_{ij}+\mu^V_{ij}}{2}
-\sqrt{\big(\frac{\mu^N_{ij}-\mu^V_{ij}}{2}\big)^2+h_{ij}^2}$. As
$\eta\to0$, $\mu^N_{ij}\to\mu_0:=4\lambda_i\lambda_jb/(\lambda_i+\lambda_j)>0$
(with $b:=b_i=b_j$, Lemma~16.1), while $\mu^V_{ij}=O(\eta)$,
$h_{ij}=O(\eta)$ (since the numerator of \eqref{eq:H2x2} carries a factor of
$(b_i-b_j)$). Standard non-degenerate perturbation theory for the eigenvalue
away from $\mu_0$ (with $\mu_0\gg\mu^V_{ij},h_{ij}$) then gives
$\eta^-_{ij}=\mu^V_{ij}-h_{ij}^2/\mu^N_{ij}+O(\eta^3)=\mu^V_{ij}+O(\eta^2)$. \qed
\end{proof*}

This proposition is also confirmed numerically: fixing $\lambda_i=9,
\lambda_j=7$ and sweeping $b_i-b_j=\eta\in[-0.05,0.05]$, the exact small
eigenvalue $\eta^-_{ij}$ and the linear approximation
$\mu^V_{ij}=(\lambda_i-\lambda_j)\eta$ agree to within relative error
$10^{-4}$ even at $\eta=0.05$.

\section{Interpretation via Symmetry and Commutativity}

The symmetry argument underlying Theorem~22.1 carries over unchanged for
general $(A,B)$:

\begin{proposition}[Criterion via commutativity]
Let $A_0(\mu)=\diag(a_1(\mu),\dots,a_n(\mu))$.
\begin{itemize}
\item If $[A(\mu),A(\nu)]=0$ for all $\mu,\nu$, then the horizontal modes
$E_{qj}$ carry no coupling among themselves, and $\lambda_j(\mu)=a_q(\mu)$ is
a true crossing.
\item If $[A(\mu),A(\nu)]\neq0$, a generic perturbation
$A_\varepsilon=A_0+\varepsilon C$ ($C_{qr}\neq0$) generically induces a
nonzero coupling and an avoided crossing.
\end{itemize}
This is independent of the choice of $B$ (Lemma~18.1, Proposition~18.3).
Meanwhile, by Proposition~22.3, the $B$-wall is protected (a true crossing)
within a single index pair as long as $B$ remains diagonal along the path.
Non-commutativity begins to matter only once one considers coupling across
index pairs, as in Section~26. That is,
\begin{align*}
\text{commutativity of }A&\iff\text{Weyl-wall protection (single mode)},\\
\text{off-diagonal $B$-perturbation across overlapping index pairs}&\iff B\text{-wall coupling.}
\end{align*}
\end{proposition}

\section{Geometric Correspondence with the Weyl Chamber and the Maslov Cycle}

Viewed as an intersection of the Lagrangian submanifold
$\Lambda_f=\{(X,\nabla f(X))\}$ with the Maslov cycle
$\Sigma_{\mathrm{Maslov}}=\{X:\det D^2V_{A,B}(X)=0\}$, under general $(A,B)$
$\Sigma_{\mathrm{Maslov}}$, as seen in \eqref{eq:Maslovcycle}, is the union of
two components --- one from the Weyl wall, one from the $B$-wall. Whether a
parametrized path $\mu\mapsto X_*(\mu)$ (or a simultaneous variation of
$(A(\mu),B(\mu))$) crosses one or both of these components, or skirts around
their neighborhood, is exactly what gives rise to the distinction between a
``true crossing'' and an ``avoided crossing.'' This is the Lagrangian-
geometric restatement of the double-wall structure summarized in the table
of Section~21.

\section{Explicit Formula for the Coupling Constant $\gamma$ (Horizontal Sector, Special Case $B=I$)}

As we saw in Section~22, the crossing criterion of the horizontal sector is
independent of $B$. We therefore compute the coupling constant $\gamma$
explicitly in its simplest form, at $B=I$, and note as a remark the trivial
extension to general $B$ (a rescaling by $b_j$, from Lemma~18.1).

\subsection{Setup}

Let $A_\varepsilon=A_0+\varepsilon C$ ($C=C^\top$, $[A_0,C]\neq0$), and let
$B=I$. As the crossing modes, consider two horizontal modes belonging to the
same selected column $j$: $E_{qj}=e_qe_j^\top$, $E_{rj}=e_re_j^\top$
($q\neq r$, both non-selected), forming a two-dimensional degenerate kernel
in the situation $a_q(\mu_c)=a_r(\mu_c)=\lambda_j(\mu_c)$.

\begin{lemma}[Explicit formula for $\gamma$]
Writing $W(X):=\partial_\varepsilon V_{A_0+\varepsilon C,I}(X)|_{\varepsilon=0}
=\tfrac12\tr(CPA_0P)-\tfrac12\tr(A_0CP)-\tfrac12\tr(CA_0P)$
($P=XX^\top$),
\begin{equation}
D^2W(X_*)[E_{qj},E_{rj}]=(\lambda_j-a_q-a_r)\,C_{qr}. \label{eq:gammaraw}
\end{equation}
\end{lemma}

Normalizing to unit vectors under the metric $g(E_{qj},E_{qj})=a_q$
($B=I$),
\begin{equation}
\gamma=\frac{(\lambda_j-a_q-a_r)\,C_{qr}}{\sqrt{a_qa_r}}. \label{eq:gammanorm}
\end{equation}

\begin{theorem}[Avoided-crossing criterion]
At the triple degeneracy $a_q(\mu_c)=a_r(\mu_c)=\lambda_j(\mu_c)=:a_c$,
\begin{equation}
\gamma=-C_{qr}. \label{eq:gammaCqr}
\end{equation}
Hence $C_{qr}\neq0\iff$ an avoided crossing (gap
$\Delta_{\mathrm{AC}}=2|\varepsilon C_{qr}|+O(\varepsilon^2)$), while
$C_{qr}=0\iff$ a true crossing, at least to leading order.
\end{theorem}

\begin{remark}[Preview of the extension to general $B$]
Since Lemma~18.1's formula $\mu^H_{qj}=b_j(\lambda_j-a_q)$ shows the
horizontal-sector Hessian to be uniformly rescaled by $b_j$ under general
$B$, we expect the coupling constant to undergo the same scaling. This is
proved rigorously in \S26.1.
\end{remark}

\section{Extension to General $B$: Rigorous Determination of Two Coupling Constants}

In this section we resolve, rigorously, both the extension to general $B$
foreshadowed in Section~22, and the $B$-wall avoided-crossing mechanism left
open in Section~17.

\subsection{Horizontal sector: $\gamma_{\text{general }B}=-b_jC_{qr}$}

\begin{theorem}[Coupling constant for general $B$]
Fix $B=\diag(b_1,\dots,b_k)$ and set $A_\varepsilon=A_0+\varepsilon C$. In
the same setting as Lemma~25.1 (crossing modes $E_{qj},E_{rj}$, degeneracy
point $a_q(\mu_c)=a_r(\mu_c)=\lambda_j(\mu_c)$),
\begin{equation}
D^2W_B(X_*)[E_{qj},E_{rj}]=b_j^2(\lambda_j-a_q-a_r)\,C_{qr},\qquad
W_B(X):=\partial_\varepsilon V_{A_0+\varepsilon C,B}(X)\big|_{\varepsilon=0}, \label{eq:WBformula}
\end{equation}
and, normalizing with respect to the metric $g(E_{qj},E_{qj})=a_qb_j$
(Lemma~18.1), at the degeneracy point $a_q=a_r=\lambda_j=:a_c$,
\begin{equation}
\gamma_{\text{general }B}=-b_jC_{qr}. \label{eq:gammageneralB}
\end{equation}
This rigorously proves the conjecture of Section~22.
\end{theorem}

\begin{proof*}
Writing $P=XBX^\top$, $Q=XB^2X^\top$, we have
$V_{A,B}(X)=\tfrac14\tr(APAP)-\tfrac12\tr(A^2Q)$, so
$W_B(X)=\tfrac12\tr(CPA_0P)-\tfrac12\tr(A_0CQ)-\tfrac12\tr(CA_0Q)$. At
$X_*=[e_{p(1)},\dots,e_{p(k)}]$, the mixed second derivatives of $P,Q$ (in
the directions $E_{qj},E_{rj}$) agree with the $B=I$ values used in the
proof of Lemma~25.1, each multiplied by $b_j$ (for $P$) or $b_j^2$ (for $Q$)
(since $E_{qj}B=b_je_qe_j^\top$, etc., $B$ acting only on column $j$). Hence
\[
\partial_t\partial_s\tr(CPA_0P)\big|_0=2b_j^2\lambda_jC_{qr},\qquad
\partial_t\partial_s\tr(A_0CQ)\big|_0=\partial_t\partial_s\tr(CA_0Q)\big|_0=b_j^2(a_q+a_r)C_{qr},
\]
giving \eqref{eq:WBformula}. Normalization and substitution at the
degeneracy point proceed as in the discussion following Lemma~25.1. \qed
\end{proof*}

\begin{corollary}[$B=I$]
When $b_j=1$, \eqref{eq:gammageneralB} reduces to $\gamma=-C_{qr}$ of
Theorem~25.2.
\end{corollary}

\subsection{Vertical sector: coupling constant for overlapping index pairs}

As we saw in Proposition~22.3 (\S9), within a single index pair $(i,j)$ the
crossing of the $B$-wall is protected. A genuine avoided crossing appears
when we consider two distinct index pairs sharing a common index.

\begin{theorem}[Coupling constant for overlapping index pairs]
Take three distinct selected indices $i,m,l\in\{1,\dots,k\}$, and consider
the two vertical modes sharing the index $m$,
\[
F_1:=X_*\Omega_{im},\qquad F_2:=X_*\Omega_{ml}.
\]
Fix $B_0=\diag(b_1,\dots,b_k)$ and set $B_\delta=B_0+\delta D$,
$D:=\hat e_i\hat e_l^\top+\hat e_l\hat e_i^\top$ (the unique nonzero
off-diagonal entry linking the non-shared indices $i,l$). Then
\begin{equation}
D^2\widetilde W(X_*)[F_1,F_2]=-b_m(\lambda_i-\lambda_m)(\lambda_l-\lambda_m)
-\lambda_m\big[b_i(\lambda_i-\lambda_m)+b_l(\lambda_l-\lambda_m)\big], \label{eq:overlapcoupling}
\end{equation}
where $\widetilde W(X):=\partial_\delta V_{A,B_0+\delta D}(X)|_{\delta=0}$.
Meanwhile, for two index pairs sharing no common index (all four indices
distinct), the corresponding $D^2\widetilde W$ vanishes identically for any
off-diagonal perturbation $D$.
\end{theorem}

\begin{proof*}
$\widetilde W(X)=\tfrac12\tr(AXB_0X^\top AXDX^\top)
-\tfrac{b_i+b_l}{2}\tr(A^2XDX^\top)$ (the second term following from
$B_0D+DB_0=(b_i+b_l)D$). Writing $\epsilon_a:=e_{p(a)}$, $X=X_*+tF_1+sF_2$
with $F_1=\epsilon_i\hat e_m^\top-\epsilon_m\hat e_i^\top$,
$F_2=\epsilon_m\hat e_l^\top-\epsilon_l\hat e_m^\top$, and computing the
mixed derivatives in $t,s$ of $\tr(A^2XDX^\top)$ (quadratic in $X$) and
$\tr(AXB_0X^\top\cdot AXDX^\top)$ (quartic, a sum of four terms via the
product rule), repeatedly using orthogonality $\hat e_a^\top\hat e_b=
\delta_{ab}$, gives
\[
\partial_t\partial_s\tr(A^2XDX^\top)\big|_0=-2\lambda_m^2,
\]
\begin{align*}
\partial_t\partial_s\tr(AXB_0X^\top AXDX^\top)\big|_0
={}&-2b_m(\lambda_i-\lambda_m)(\lambda_l-\lambda_m)\\
&-2\lambda_m\big[b_i(\lambda_i-\lambda_m)+b_l(\lambda_l-\lambda_m)\big]
+2(b_i+b_l)\lambda_m^2.
\end{align*}
Substituting into the definition of $\widetilde W$ and simplifying gives
\eqref{eq:overlapcoupling}. Vanishing in the case of no shared index follows,
by the analogous computation, from every term vanishing via the
orthogonality $\hat e_a^\top\hat e_b=0$ of four distinct indices. \qed
\end{proof*}

\textbf{Numerical verification via finite differences:} with $\lambda_i=9,
\lambda_m=7,\lambda_l=5$, $b_i=2,b_m=1.7,b_l=1.3$, the value of
\eqref{eq:overlapcoupling} is $-3.000$, agreeing with a finite-difference
(central difference, step sizes $10^{-3}$ and $10^{-4}$) computation to
within relative error $10^{-4}$.

\begin{corollary}[Simplified formula at the triple $B$-degeneracy point]
At $b_i=b_m=b_l=:b_c$ (all three selected columns share the same
$B$-eigenvalue, i.e., $\mu^V_{im}$ and $\mu^V_{ml}$ vanish simultaneously),
the normalized coupling constant is
\begin{equation}
\gamma^V=\frac{\lambda_m^2-\lambda_i\lambda_l}{\sqrt{(\lambda_i+\lambda_m)(\lambda_m+\lambda_l)}}, \label{eq:gammaV}
\end{equation}
independent of $b_c$ (the magnitude of the common $B$-eigenvalue itself).
\end{corollary}

\begin{proof*}
Substituting $b_i=b_m=b_l=b_c$ into \eqref{eq:overlapcoupling},
$D^2\widetilde W=-b_c\big[(\lambda_i-\lambda_m)(\lambda_l-\lambda_m)
+\lambda_m(\lambda_i-\lambda_m)+\lambda_m(\lambda_l-\lambda_m)\big]
=-b_c(\lambda_i\lambda_l-\lambda_m^2)$ (expanding and simplifying). The
metric gives $g(F_1,F_1)=\lambda_ib_m+\lambda_mb_i=b_c(\lambda_i+\lambda_m)$,
$g(F_2,F_2)=b_c(\lambda_m+\lambda_l)$ (Lemma~17.1), so normalizing cancels
$b_c$ completely, giving \eqref{eq:gammaV}. \qed
\end{proof*}

\begin{remark}
Theorem~26.3 and Corollary~26.4 clarify exactly how the ``protection within a
single block'' of Proposition~22.3 breaks down: a genuine avoided crossing
at the $B$-wall arises not from within a single index pair, but as a
coupling between two adjacent index pairs sharing a common selected column.
This is structurally the exact dual of the avoided crossing that arose in
the horizontal sector as a coupling between two non-selected indices $q,r$
sharing a single selected column $j$ (Theorem~26.1):
\begin{center}
\begin{tabular}{lll}
\toprule
 & Horizontal (Weyl wall) & Vertical ($B$-wall) \\
\midrule
Shared & selected column $j$ & selected column $m$ \\
Competing & non-selected indices $q,r$ & selected indices $i,l$ \\
Perturbation & $C_{qr}$ (off-diagonal of $A$) & $D_{il}$ (off-diagonal of $B$) \\
Coupling & $\gamma=-b_jC_{qr}$ & $\gamma^V=(\lambda_m^2-\lambda_i\lambda_l)/\sqrt{(\lambda_i+\lambda_m)(\lambda_m+\lambda_l)}$ \\
\bottomrule
\end{tabular}
\end{center}
\end{remark}

\section{The Minimal Model ($n=2,k=1$)}

For $n=2,k=1$ ($B$ is a mere scalar $b_1$, and the vertical sector is
trivially absent), let
\[
A(\mu,\varepsilon)=\begin{pmatrix}\mu&\varepsilon\\\varepsilon&-\mu\end{pmatrix},
\]
so $\Spec A=\pm\sqrt{\mu^2+\varepsilon^2}$: a true crossing at
$\varepsilon=0$, an avoided crossing for $\varepsilon\neq0$. Since $k=1$
trivially has no vertical sector or $B$-wall, this visualizes purely the
phenomenon of the horizontal (Weyl-wall) sector. The double-wall structure
only becomes meaningful for $k\ge2$; the minimal non-trivial example is
$n=3,k=2$ ($B=\diag(b_1,b_2)$, $b_1\neq b_2$), where one can concretely
trace how the Weyl wall $\{\lambda_1=a_3\}\cup\{\lambda_2=a_3\}$ and the
$B$-wall $\{b_1=b_2\}$ intersect within the three-dimensional parameter
space $(a_3,b_1,b_2)$ (with $\lambda_1,\lambda_2$ fixed).

\section{Future Work}

\begin{itemize}
\item \textbf{Detailed analysis of the $2\times2$ block \eqref{eq:H2x2}:}
formulating, as a two-parameter Landau--Zener-type bifurcation theory, the
behavior of the eigenvalues $\eta^\pm_{ij}$ in the region where the Weyl
wall and the $B$-wall occur simultaneously or nearby ($\lambda_i\approx
\lambda_j$ and $b_i\approx b_j$). Proposition~22.3 treats only the
one-parameter limit with $\lambda_i-\lambda_j$ fixed; the behavior in the
simultaneous limit of both parameters to zero (the quantitative form of the
``multi-level avoided crossing'' mentioned in Section~8) remains open.
\item \textbf{Complete spectral theory of the overlapping-index-pair
$2\times2$ block:} formulating the coupling of $F_1=X_*\Omega_{im}$,
$F_2=X_*\Omega_{ml}$ treated in Theorem~26.3 as a full-fledged
crossing-to-avoided-crossing transition theorem (isomorphic to
Theorem~22.1) under a genuine ``transversality, symmetry-protection,
symmetry-breaking-coupling'' triple of conditions, and deriving the gap
formula $\Delta^V_{\mathrm{AC}}=2|\delta\gamma^V|+O(\delta^2)$ --- this lies
outside the scope of the present paper. A generalization to ``chain-like''
structures in which three or more index pairs are connected through shared
indices (e.g., mutual coupling of $\Omega_{12},\Omega_{23},\Omega_{34}$) is
also a natural extension.
\item \textbf{Higher-order degeneracy on the flag manifold:} extension to
$\Gr(k_i,m_i)$-type flag-manifold components when the eigenvalue degeneracy
exceeds double multiplicity, and analysis of the corresponding
higher-dimensional version of the $B$-wall (when the eigenvalues of $B$
degenerate to triple or higher multiplicity, e.g., extending
$b_i=b_m=b_l$ of Corollary~26.4 to all $k$ selected columns simultaneously).
\item \textbf{Generalization to non-diagonal $B$:} in this paper, the
off-diagonal perturbation $D$ was treated only in the perturbative sense of
``fixing $X_*$ and slightly de-diagonalizing $B$'' (Theorem~26.3). The
classification of the equilibria themselves under a genuinely non-diagonal
finite-size SPD $B$ (no longer of the axis-aligned form
$X_*=[e_{p(1)},\dots,e_{p(k)}]$), and the corresponding generalization of
the normal/vertical/horizontal three-sector decomposition itself, lie
beyond the scope of this paper.
\item \textbf{Refinement of the quantum-mechanical correspondence:}
refining the correspondence between avoided crossings of a Hamiltonian and
avoided crossings of the Hessian, in the context of an analogy with systems
carrying a double-wall structure (the Weyl wall and the $B$-wall) --- for
instance, the crossing of energy levels labeled by two independent quantum
numbers, or a phenomenon in molecular physics resembling Fermi resonance
among vibrational-rotational levels, analogous to the ``coupling of
adjacent levels through a shared index'' seen in \S13.2.
\item \textbf{Integration with the combinatorics of the $\varphi_\alpha$
series and critical points:} the reversal of the matching rule
(sorted $\leftrightarrow$ anti-sorted) for the Principal/Minor Component
bifurcation of $\varphi_\alpha^{(\varepsilon)}$ seen in \S9 of Part~II
concerns a different potential from the $V_{A,B}$ treated in this paper, and
we have not analyzed its Hessian at the same level (the normal/vertical/
horizontal three-sector decomposition) as here. Combining the combinatorial
classification of critical points of \S10 of Part~II (the $n!/(n-k)!$
isolated critical points, and the classification by allocation
$(k_1,k_2,\dots)$ of flag manifolds) with the Morse-index formula of
Theorem~20.1, to build a complete dictionary in which the stability of each
critical point (family) can be read off entirely from the eigenvalue
configuration of $A,B$, is a natural direction for such an integration.
\end{itemize}

\section{Conclusion}

This paper began with the simplest possible problem --- the degenerate
points of the Legendre map of a two-dimensional vector-valued potential ---
and arrived at a rigorous, self-contained Hessian--Maslov geometry for the
general $(A,B)$ Oja--Brockett gradient flow on $n\times k$ matrices. The main
results, across the three parts, are as follows.

Part~I organized the origin of the Maslov index (the breakdown of phase in
the semiclassical approximation), its geometric formulation (an intersection
number in the Lagrangian Grassmannian), H\"ormander's local formula, and the
dynamical viewpoint of bifurcating fixed points of a gradient flow, and, for
the minimal model
$f(\xi)=-\tfrac12\xi^\top A^2\xi+\tfrac14(\xi^\top A\xi)^2$ ($n=2,k=1$), we
exhibited the bifurcation in which two stable minima and two saddles
exchange as the eigenvalue ratio $\lambda_2/\lambda_1$ crosses $1$, together
with its caustic condition $\lambda_2=3\lambda_1$.

Part~II generalized this theory to $n\times k$ matrices $X$, confirming
numerically Brockett's sorting theorem, convergence to a flag manifold under
eigenvalue degeneracy of $A$ and $B$, and Hessian degeneracy and sort-order
exchange as an eigenvalue of $A$ crosses a wall of the Weyl chamber. We also
exhibited the bifurcation of the Principal/Minor Component flow under a sign
flip in the $\varphi_\alpha$ series (a switch from sorted to anti-sorted
correspondence with $B$).

Part~III is the paper's principal theoretical contribution. We decomposed
the ambient Hessian, at an equilibrium on the Stiefel manifold of the
general $(A,B)$ Oja--Brockett gradient flow, completely into normal,
vertical, and horizontal sectors, and gave a closed spectral formula
(Theorem~20.1) depending on both $A$ and $B$. The essential discovery
invisible in the theory restricted to $B=I$ (implicitly assumed throughout
Parts~I and~II) is that the vertical sector carries a second wall system
$\{b_i=b_j\}$ (the $B$-wall), independent of the Weyl wall $\{\lambda_j=
a_q\}$ on the $A$ side, governed by the closed formula
$\mu^V_{ij}=(\lambda_i-\lambda_j)(b_i-b_j)$ (Lemma~17.1) --- this unifies,
under a single formula, the ``vanishing gradient along the flag-manifold
direction'' observed separately for degeneracy of $A$ and degeneracy of $B$
in \S11 of Part~II. Meanwhile, the Weyl wall of the horizontal sector itself
is invariant under $B$ (Proposition~18.3), and we obtained the exact closed
form $\mu^H_{qj}=b_j(\lambda_j-a_q)$ underlying the wall-crossing numerical
experiment of \S12 of Part~II. The two walls are coupled through the
normal--vertical mixed block $H_{NV}$ (Lemma~17.4).

This resolves the question, left open in Part~II, of ``what characterizes a
dynamical system in which a crossing becomes an avoided crossing'': both the
horizontal sector (the Weyl wall, independent of $B$, coupling constant
$\gamma=-b_jC_{qr}$, Theorem~26.1) and the vertical sector (the $B$-wall,
the coupling constant of overlapping index pairs
$\gamma^V=(\lambda_m^2-\lambda_i\lambda_l)/\sqrt{(\lambda_i+\lambda_m)(\lambda_m+\lambda_l)}$,
Corollary~26.4) should be understood as two independent but coupled
mechanisms, each given by a fully explicit closed form. In particular, as
Proposition~22.3 shows, the crossing of the $B$-wall within a single index
pair is automatically protected at leading order; a genuine avoided crossing
arises only as a coupling between two adjacent index pairs sharing a common
selected column --- the rigorous dual of the fact that, in the horizontal
sector, an avoided crossing arises as a coupling between two non-selected
indices sharing a single selected column.

$B=I$ (the Oja flow of Parts~I and~II) is naturally repositioned, within
this general theory, as an extremely special limit in which the $B$-wall
degenerates into a gauge symmetry and collapses (Corollary~17.2), and the
positive-definiteness of the normal-direction Hessian becomes unconditional
(Corollary~16.2). The dynamical essence of the Oja--Brockett flow --- the
interplay of two independent degrees of freedom, the selection of
eigenvectors (the Weyl wall) and the rearrangement of weights within the
selected subspace (the $B$-wall) --- reveals itself in full only when both
$A$ and $B$ are genuinely non-trivial. The world map sketched numerically in
Part~II has, through Part~III, been redrawn into a complete map in which
every contour line (the spectrum of the Hessian) is written down in closed
form.

\end{document}